\documentclass[12pt]{article}
\usepackage{amsmath,amssymb}
\usepackage{graphicx}
\newcommand{\eqdef}{\stackrel{\text{def}}{=}}
\newcommand{\n}{\nonumber\\}
\newcommand{\bm}{\boldsymbol}
\newcommand{\ignore}[1]{}
\numberwithin{equation}{section}
\newcommand{\Romannumeral}[1]{\uppercase\expandafter{\romannumeral#1}}
\newcommand{\I}{\text{\Romannumeral{1}}}
\newcommand{\II}{\text{\Romannumeral{2}}}
\newcommand{\III}{\text{\Romannumeral{3}}}
\newcommand{\norm}[1]{|\!|#1|\!|}

\allowdisplaybreaks[4]

\begin{document}

\baselineskip=20pt

\newcommand{\preprint}{
\vspace*{-20mm}
   \begin{flushright}\normalsize \sf
    DPSU-16-3\\
  \end{flushright}}
\newcommand{\Title}[1]{{\baselineskip=26pt
  \begin{center} \Large \bf #1 \\ \ \\ \end{center}}}
\newcommand{\Author}{\begin{center}
  \large \bf Satoru Odake and Ryu Sasaki \end{center}}
\newcommand{\Address}{\begin{center}
     Faculty of Science, Shinshu University,\\
     Matsumoto 390-8621, Japan
   \end{center}}
\newcommand{\Accepted}[1]{\begin{center}
  {\large \sf #1}\\ \vspace{1mm}{\small \sf Accepted for Publication}
  \end{center}}

\preprint
\thispagestyle{empty}

\Title{Multi-indexed Meixner and Little $q$-Jacobi (Laguerre) Polynomials}

\Author

\Address
\vspace{1cm}

\begin{abstract}
As the fourth stage of the project {\em multi-indexed orthogonal
polynomials\/}, we present the multi-indexed Meixner and little $q$-Jacobi
(Laguerre) polynomials in the framework of `discrete quantum mechanics' with
real shifts defined on the semi-infinite lattice in one dimension.
They are obtained, in a similar way to the multi-indexed Laguerre and Jacobi
polynomials reported earlier, from the quantum mechanical systems
corresponding to the original orthogonal polynomials by multiple application
of the discrete analogue of the Darboux transformations or the
Crum-Krein-Adler deletion of virtual state vectors.
The virtual state vectors are the solutions of the matrix Schr\"odinger
equation on all the lattice points having negative energies and infinite norm.
This is in good contrast to the ($q$-)Racah systems defined on a finite
lattice, in which the `virtual state' vectors satisfy the matrix
Schr\"odinger equation except for one of the two boundary points.
\end{abstract}

\section{Introduction}
\label{sec:intro}

Theory of exactly solvable quantum mechanics (QM) has seen a great surge of
interests in recent years aroused by the discovery of the {\em exceptional\/}
\cite{gomez,quesne,os16} and {\em multi-indexed orthogonal polynomials\/}
\cite{gomez2,os25}. They are new types of orthogonal polynomials satisfying
second order differential equations and forming complete sets of orthogonal
bases in appropriate Hilbert spaces. They are distinguished from the
{\em classical orthogonal polynomials} \cite{askey}, by the fact that there
are `holes' in their degrees and the breakdown of the three term recurrence
relations, which is essential to avoid the so-called
Bochner's theorem \cite{routh,bochner}.

The concept of the exceptional and multi-indexed orthogonal polynomials has
been generalised to the other classical orthogonal polynomials in the Askey
scheme \cite{kls}, which satisfy second order
{\em difference\/} equations. The exceptional and multi-indexed Wilson,
Askey-Wilson, ($q$-)Racah polynomials have been constructed \cite{os27,os26}
as the main parts of the eigenfunctions (vectors) of difference Schr\"odinger
equations within the program of `{\em discrete quantum mechanics\/}'
\cite{os24,os12,os13}.

In this paper we report the construction of the multi-indexed Meixner,
little $q$-Jacobi (Laguerre) polynomials.
The base polynomials, the Meixner (M), little $q$-Jacobi (l$q$J) and little
$q$-Laguerre (l$q$L) polynomials, satisfy second order difference equations
on {\em semi-infinite\/} integer lattices and they form complete orthogonal
bases in the corresponding $\ell^2$ Hilbert spaces.
These polynomials are also called `orthogonal polynomials of a discrete
variable' \cite{nikiforov} and they belong to the subclass of discrete
quantum mechanics with real shifts.
The Meixner case was studied by Dur\'an \cite{d14}.
In our language his polynomials correspond to the eigenstates and/or virtual
states deletion.

The first stage of the multi-indexed orthogonal polynomials project dealt
with the Laguerre and the Jacobi polynomials \cite{os25} belonging to the
ordinary QM.
The multi-indexed Wilson and Askey-Wilson polynomials were constructed in
the second stage \cite{os27}. They belong to discrete QM with pure imaginary
shifts \cite{os24,os13}.
The multi-indexed ($q$-)Racah polynomials were constructed in the third stage
\cite{os26}. They belong to the discrete QM with real shifts \cite{os12,os34}
defined on {\em finite\/} integer lattices.
The multi-indexed orthogonal polynomials are obtained as the main parts of
the eigenfunctions of the deformed quantum systems, which are generated by
multiple application of Darboux-Crum transformations \cite{darboux,crum} on
the base systems.
The Darboux-Crum transformations for discrete QM with real shifts were
developed in \cite{os22}.
The seed solutions for the Darboux-Crum transformations
are called the {\em virtual state functions (vectors)\/}, which are obtained
from the eigenfunctions by the discrete symmetry transformations of the
base Hamiltonians.

The present paper is organised as follows. In section two, the general
setting of the discrete QM with real shifts is briefly recapitulated.
The main features of the base systems, the Meixner, little $q$-Jacobi
(Laguerre) polynomials are collected in \S\ref{sec:org}.
The method of constructing virtual state vectors is explained in
\S\ref{sec:virt} together with their explicit forms.
The general procedure of modifying exactly solvable QM by multiple
Darboux-Crum transformations in terms of the virtual state vectors is
explained in some detail in section \ref{sec:mvsd}.
Various formulas of the obtained multi-indexed M, l$q$J, l$q$L polynomials
are collected in section \ref{sec:miop}.
Certain limits of these multi-indexed polynomials are presented in
\S\,\ref{sec:limit}.
The final section is for a summary and comments.

\goodbreak
\section{Discrete QM with Real Shifts on a Semi-infinite Lattice}
\label{sec:rdQM}

Let us recapitulate the discrete quantum mechanics with real shifts on a
semi-infinite lattice.
There are two types:
one component systems (Meixner, little $q$-Jacobi, etc.) \cite{os12}
and two component systems ($q$-Meixner, big $q$-Jacobi, etc.) \cite{os34}.
We discuss one component systems in this paper.

The Hamiltonian $\mathcal{H}=(\mathcal{H}_{x,y})$ is a {\em tri-diagonal real
symmetric\/} (Jacobi) matrix and its rows and columns are indexed by
integers $x$ and $y$, which take values in $\mathbb{Z}_{\geq 0}$
(semi-infinite).
By adding a scalar matrix to the Hamiltonian, the lowest eigenvalue is
adjusted to be zero. This makes the Hamiltonian {\em positive semi-definite}.
Since the eigenvector corresponding to the zero eigenvalue has definite sign,
{\em i.e.\/} all the components are positive or negative, the Hamiltonian
$\mathcal{H}$ has the following form ($x,y\in\mathbb{Z}_{\ge0}$)
\begin{equation}
  \mathcal{H}_{x,y}\eqdef
  -\sqrt{B(x)D(x+1)}\,\delta_{x+1,y}-\sqrt{B(x-1)D(x)}\,\delta_{x-1,y}
  +\bigl(B(x)+D(x)\bigr)\delta_{x,y},
  \label{Hdef}
\end{equation}
in which the potential functions $B(x)$ and $D(x)$ are real and positive
but vanish at the boundary
\begin{equation}
  B(x)>0\ \ (x\in\mathbb{Z}_{\geq 0}),\quad
  D(x)>0\ \ (x\in\mathbb{Z}_{\geq 1}),\ \ D(0)=0.
  \label{BDcond}
\end{equation}
The potential functions $B(x)$ and $D(x)$ are rational functions of $x$ or
$q^x$ ($q$ is a positive constant, $0<q<1$).
See \S\ref{sec:org} for the explicit forms of these functions.
The Hamiltonian \eqref{Hdef} is real symmetric, $\mathcal{H}_{x,y}
=\mathcal{H}_{y,x}$.
Reflecting the positive semi-definiteness, the Hamiltonian \eqref{Hdef}
can be expressed in a factorised form:
\begin{align}
  &\mathcal{H}=\mathcal{A}^{\dagger}\mathcal{A},\qquad
  \mathcal{A}=(\mathcal{A}_{x,y}),
  \ \ \mathcal{A}^{\dagger}=((\mathcal{A}^{\dagger})_{x,y})
  =(\mathcal{A}_{y,x}),
  \label{factor}\\
  &\mathcal{A}_{x,y}\eqdef
  \sqrt{B(x)}\,\delta_{x,y}-\sqrt{D(x+1)}\,\delta_{x+1,y},\quad
  (\mathcal{A}^{\dagger})_{x,y}=
  \sqrt{B(x)}\,\delta_{x,y}-\sqrt{D(x)}\,\delta_{x-1,y}.
\end{align}
Here $\mathcal{A}$ ($\mathcal{A}^\dagger$) is an upper (lower) triangular
matrix with the diagonal and the super(sub)-diagonal entries only.
For simplicity in notation, we write $\mathcal{H}$, $\mathcal{A}$ and
$\mathcal{A}^{\dagger}$ as follows:
\begin{align}
  &e^{\pm\partial}=((e^{\pm\partial})_{x,y}),\quad
  (e^{\pm\partial})_{x,y}\eqdef\delta_{x\pm 1,y},\quad
  (e^{\partial})^{\dagger}=e^{-\partial},
  \label{partdef}\\
  &\mathcal{H}=-\sqrt{B(x)D(x+1)}\,e^{\partial}
  -\sqrt{B(x-1)D(x)}\,e^{-\partial}+B(x)+D(x)\n
  &\phantom{\mathcal{H}}=-\sqrt{B(x)}\,e^{\partial}\sqrt{D(x)}
  -\sqrt{D(x)}\,e^{-\partial}\sqrt{B(x)}+B(x)+D(x),
  \label{Hdef2}\\
  &\mathcal{A}=\sqrt{B(x)}-e^{\partial}\sqrt{D(x)},\quad
  \mathcal{A}^{\dagger}=\sqrt{B(x)}-\sqrt{D(x)}\,e^{-\partial}.
  \label{A,Ad}
\end{align}
We suppress the unit matrix $\bm{1}=(\delta_{x,y})$:
$(B(x)+D(x))\bm{1}$ in \eqref{Hdef2}, $\sqrt{B(x)}\,\bm{1}$ in \eqref{A,Ad}.
Note that the product of $e^{\partial}$ and $e^{-\partial}$ is
$e^{\partial}e^{-\partial}=\bm{1}$ but the reversed order is
$e^{-\partial}e^{\partial}=\bm{1}-\text{diag}(1,0,0,\ldots)$.

The Hamiltonian \eqref{Hdef} is a linear operator on the real $\ell^2$
Hilbert space with the inner product of two real vectors $f=(f(x))$ and
$g=(g(x))$ defined by
\begin{equation}
  (f,g)\eqdef\lim_{N\to\infty}(f,g)_N,\quad
  (f,g)_N\eqdef\sum_{x=0}^Nf(x)g(x),\quad
  \norm{f}^2\eqdef(f,f)<\infty.
  \label{inpro}
\end{equation}
The Schr\"odinger equation is the eigenvalue problem for the hermitian
matrix $\mathcal{H}$,
\begin{equation}
  \mathcal{H}\phi_n(x)=\mathcal{E}_n\phi_n(x)\quad
  (n=0,1,\ldots),\quad
  0=\mathcal{E}_0<\mathcal{E}_1<\cdots,
  \label{schreq}
\end{equation}
where the eigenvector $\phi_n=(\phi_n(x))$ is, by definition, of finite
norm, $\norm{\phi_n}<\infty$.
Let us recall the fact that the spectrum of a Jacobi matrix is simple.
The ground state eigenvector, which is chosen positive $\phi_0(x)>0$
($x\in\mathbb{Z}_{\geq 0}$), satisfies the zero mode equation:
\begin{equation}
  \mathcal{A}\phi_0=0\ \Rightarrow\mathcal{H}\phi_0=0,\quad
  \sqrt{B(x)}\,\phi_0(x)=\sqrt{D(x+1)}\,\phi_0(x+1),
  \label{phi0eq}
\end{equation}
and it is easily obtained with the normalisation $\phi_0(0)=1$
(convention: $\prod\limits_{k=n}^{n-1}*=1$):
\begin{equation}
  \phi_0(x)=\prod_{y=0}^{x-1}\sqrt{\frac{B(y)}{D(y+1)}}\,.
  \label{phi0}
\end{equation}
The self-adjointness of the Hamiltonian and the non-degeneracy of the
spectrum \eqref{schreq} imply that the eigenvectors are mutually orthogonal:
\begin{equation}
  (\phi_n,\phi_m)=\frac{\delta_{nm}}{d_n^2}\quad
  (n,m\in\mathbb{Z}_{\geq 0}),
  \label{ortho}
\end{equation}
where $d_n^2$ ($d_n>0$) is the normalisation constant.

For the systems considered in this paper, the eigenvectors have the following
factorised form
\begin{equation}
  \phi_n(x)=\phi_0(x)\check{P}_n(x),\quad
  \check{P}_n(x)\eqdef P_n\bigl(\eta(x)\bigr),
\end{equation}
where $\eta(x)$ is the {\em sinusoidal coordinate} satisfying the boundary
condition $\eta(0)=0$.
The other function $P_n(\eta)$ is a degree $n$ polynomial in $\eta$ and
we adopt the universal normalisation condition \cite{os12,os34} as
\begin{equation}
  P_n(0)=1\ \bigl(\Leftrightarrow \check{P}_n(0)=1\bigr).
  \label{Pzero}
\end{equation}
This is important for the proper definition of the {\em dual polynomials\/},
in which the roles of $x$ and $n$ are interchanged \cite{os12}.
The sinusoidal coordinate has a special dynamical meaning \cite{os7,os12}.
The Heisenberg operator solution for $\eta(x)$ can be expressed
in a closed form, which is a consequence of the closure relation \cite{os12}.
This means that its time evolution is a sinusoidal motion.
The similarity transformed Hamiltonian $\widetilde{\mathcal{H}}$ in terms of
the ground state wavefunction $\phi_0(x)$ \eqref{phi0} is
\begin{equation}
  \widetilde{\mathcal{H}}
  \eqdef\phi_0(x)^{-1}\circ\mathcal{H}\circ\phi_0(x)
  =B(x)(1-e^{\partial})+D(x)(1-e^{-\partial}),
  \label{Ht}
\end{equation}
and \eqref{schreq} becomes square root free
\begin{equation}
  \widetilde{\mathcal{H}}\check{P}_n(x)=\mathcal{E}_n\check{P}_n(x),
  \label{Hteq}
\end{equation}
namely
\begin{equation}
  B(x)\bigl(\check{P}_n(x)-\check{P}_n(x+1)\bigr)
  +D(x)\bigl(\check{P}_n(x)-\check{P}_n(x-1)\bigr)
  =\mathcal{E}_n\check{P}_n(x),
  \label{tHcPn=}
\end{equation}
which are valid for any $x$($\in\mathbb{C}$), as $B(x)$ and $D(x)$ are
rational functions of $x$ or $q^x$.

\subsection{Original systems}
\label{sec:org}

Let us consider the Meixner (M), little $q$-Jacobi (l$q$J) and little
$q$-Laguerre (l$q$L) cases.
We follow the notation of \cite{os12}.
Various quantities depend on a set of parameters
$\bm{\lambda}=(\lambda_1,\lambda_2,\ldots)$ and their dependence is
expressed like,
$\mathcal{H}=\mathcal{H}(\bm{\lambda})$,
$\mathcal{A}=\mathcal{A}(\bm{\lambda})$,
$\mathcal{E}_n=\mathcal{E}_n(\bm{\lambda})$,
$B(x)=B(x;\bm{\lambda})$,
$\phi_n(x)=\phi_n(x;\bm{\lambda})$, 
$\check{P}_n(x)=\check{P}_n(x;\bm{\lambda})$, etc.
The parameter $q$ is $0<q<1$ and $q^{\bm{\lambda}}$ stands for
$q^{(\lambda_1,\lambda_2,\ldots)}=(q^{\lambda_1},q^{\lambda_2},\ldots)$.
For these three systems the sinusoidal coordinate $\eta(x)$, the auxiliary
function $\varphi(x)$ \eqref{Meixnereta}, \eqref{lqJeta} and \eqref{lqLeta}
and the potential function $D(x)$ do not depend on the parameters
$\bm{\lambda}$.
The normalisation (boundary condition) of $\eta(x)$ is chosen as $\eta(0)=0$.
The shift parameter $\bm{\delta}$ and the scale parameter $\kappa$ appear
in the defining formula of shape invariance \eqref{shinv}.

\subsubsection{Meixner (M)}
\label{sec:M}

We rescale the overall normalisation of the Hamiltonian in \cite{os12}:
$(\mathcal{H}\text{ in \cite{os12}})\times(1-c)\to\mathcal{H}$.
The fundamental data are as follows \cite{os12}:
\begin{align}
  &\bm{\lambda}=(\beta,c),\quad
  \bm{\delta}=(1,0),\quad \kappa=1,\quad \beta>0,\quad 0<c<1,\\
  &B(x;\bm{\lambda})=c(x+\beta),\quad
  D(x)=x,
  \label{MeixnerBD}\\
  &\mathcal{E}_n(\bm{\lambda})=(1-c)n,\quad\eta(x)=x,\quad
  \varphi(x)=1,
  \label{Meixnereta}\\
  &\check{P}_n(x;\bm{\lambda})
  ={}_2F_1\Bigl(
  \genfrac{}{}{0pt}{}{-n,\,-x}{\beta}\Bigm|1-c^{-1}\Bigr)
  =M_n(x;\beta,c),
  \label{MeixnerP}\\
  &\phi_0(x;\bm{\lambda})^2=\frac{(\beta)_x\,c^x}{x!},\quad
  d_n(\bm{\lambda})^2
  =\frac{(\beta)_n\,c^n}{n!}\times(1-c)^{\beta},
  \label{Meixnerphi0d}
\end{align}
where $M_n(\eta;\beta,c)$ is the Meixner polynomial,
which is obviously {\em self-dual\/}, $M_n(x;\beta,c)=M_x(n;\beta,c)$
$(x,n\in\mathbb{Z}_{\geq 0})$.

\subsubsection{little $q$-Jacobi (l$q$J)}
\label{sec:lqJ}

The little $q$-Jacobi polynomial to be discussed in this paper \eqref{lqJuni}
obeys the universal normalisation \eqref{Pzero}, which is different from the
conventional normalisation of the little $q$-Jacobi polynomial $p_n$ as shown
explicitly in \eqref{littleqjacobinorm}.
The fundamental data are as follows \cite{os12,os34}:
\begin{align}
  &q^{\bm{\lambda}}=(a,b),\quad
  \bm{\delta}=(1,1),\quad \kappa=q^{-1},\quad 0<a<q^{-1},\quad b<q^{-1},\\
  &B(x;\bm{\lambda})=a(q^{-x}-bq),\quad
  D(x)=q^{-x}-1,\\
  &\mathcal{E}_n(\bm{\lambda})=(q^{-n}-1)(1-abq^{n+1}),\quad
  \eta(x)=1-q^x,\quad\varphi(x)=q^x,
  \label{lqJeta}\\
  &\check{P}_n(x;\bm{\lambda})
  ={}_3\phi_1\Bigl(
  \genfrac{}{}{0pt}{}{q^{-n},\,abq^{n+1},\,q^{-x}}{bq}\Bigm|
  q\,;\frac{q^x}{a}\Bigr)
  \label{lqJuni}\\
  &\phantom{\check{P}_n(x;\bm{\lambda})}
  =\frac{(a^{-1}q^{-n};q)_n}{(bq;q)_n}\,
  {}_2\phi_1\Bigl(
  \genfrac{}{}{0pt}{}{q^{-n},\,abq^{n+1}}{aq}\Bigm|q\,;q^{x+1}\Bigr)\n
  &\phantom{\check{P}_n(x;\bm{\lambda})}
  =\frac{(a^{-1}q^{-n};q)_n}{(bq;q)_n}\,
  p_n\bigl(1-\eta(x);a,b|q\bigr),
  \label{littleqjacobinorm}\\
  &\phi_0(x;\bm{\lambda})^2=\frac{(bq;q)_x}{(q;q)_x}(aq)^x,
  \ \ d_n(\bm{\lambda})^2
  =\frac{(bq,abq;q)_n\,a^nq^{n^2}}{(q,aq;q)_n}\,
  \frac{1-abq^{2n+1}}{1-abq}
  \times\frac{(aq;q)_{\infty}}{(abq^2;q)_{\infty}}.
  \label{littleqjacobidn}
\end{align}
For the dual l$q$J polynomial, see \cite{atakishi2} and \cite{os12}.

\subsubsection{little $q$-Laguerre (l$q$L)}
\label{sec:lqL}

The little $q$-Laguerre polynomial to be discussed in this paper
\eqref{lqLuni} obeys the universal normalisation \eqref{Pzero}, which is
different from the conventional normalisation of the little $q$-Laguerre
polynomial $p_n$ as shown explicitly in \eqref{littleqlaguerrenorm}.
The fundamental data are as follows \cite{os12}:
\begin{align}
  &q^{\bm{\lambda}}=a,\quad
  \bm{\delta}=1,\quad \kappa=q^{-1},\quad 0<a<q^{-1},\\
  &B(x;\bm{\lambda})=aq^{-x},\quad
  D(x)=q^{-x}-1,\\
  &\mathcal{E}_n=q^{-n}-1,\quad
  \eta(x)=1-q^x,\quad\varphi(x)=q^x,
  \label{lqLeta}\\
  &\check{P}_n(x;\bm{\lambda})
  ={}_2\phi_0\Bigl(
  \genfrac{}{}{0pt}{}{q^{-n},\,q^{-x}}{-}\Bigm|q;\frac{q^x}{a}\Bigr)
  =(a^{-1}q^{-n}\,;q)_n\,{}_2\phi_1\Bigl(
  \genfrac{}{}{0pt}{}{q^{-n},\,0}{aq}\Bigm|q;q^{x+1}\Bigr)
  \label{lqLuni}\\
  &\phantom{\check{P}_n(x;\bm{\lambda})}
  =(a^{-1}q^{-n}\,;q)_n\,p_n\bigl(1-\eta(x);a|q\bigr),
  \label{littleqlaguerrenorm}\\
  &\phi_0(x;\bm{\lambda})^2=\frac{(aq)^x}{(q;q)_x},\quad
  d_n(\bm{\lambda})^2
  =\frac{a^nq^{n^2}}{(q,aq;q)_n}\times(aq;q)_{\infty}.
  \label{littleqlaguerredn}
\end{align}
The l$q$L system is obtained from l$q$J system by setting $b=0$.
The dual l$q$L polynomial is the Al-Salam-Carlitz $\II$ polynomial
\cite{kls,atakishi2,os12}.

\subsubsection{shape invariance}
\label{sec:si}

These three systems are shape invariant \cite{os12},
\begin{equation}
  \mathcal{A}(\bm{\lambda})\mathcal{A}(\bm{\lambda})^{\dagger}
  =\kappa\mathcal{A}(\bm{\lambda}+\bm{\delta})^{\dagger}
  \mathcal{A}(\bm{\lambda}+\bm{\delta})+\mathcal{E}_1(\bm{\lambda}),
  \label{shinv}
\end{equation}
in which $\bm{\delta}$ denotes the shift of the parameters, $\kappa$ is
a positive constant and $\mathcal{E}_1(\bm{\lambda})$ is the eigenvalue
of the first excited state $\mathcal{E}_1(\bm{\lambda})>0$.
It is a sufficient condition for exact solvability and it provides the
explicit formulas for the energy eigenvalues
$\mathcal{E}_n(\bm{\lambda})=\sum\limits_{s=0}^{n-1}\kappa^s
\mathcal{E}_1(\bm{\lambda}+s\bm{\delta})$ and the eigenfunctions,
{\em i.e.} the universal Rodrigues formula \cite{os12,os34}.
The forward and backward shift relations are
\begin{equation}
  \mathcal{F}(\bm{\lambda})\check{P}_n(x;\bm{\lambda})
  =\mathcal{E}_n(\bm{\lambda})\check{P}_{n-1}(x;\bm{\lambda}+\bm{\delta}),
  \quad
  \mathcal{B}(\bm{\lambda})\check{P}_{n-1}(x;\bm{\lambda}+\bm{\delta})
  =\check{P}_n(x;\bm{\lambda}),
  \label{BPn-1=Pn}
\end{equation}
where the forward and backward shift operators are
\begin{align}
  \mathcal{F}(\bm{\lambda})
  &\eqdef\sqrt{B(0;\bm{\lambda})}\,\phi_0(x;\bm{\lambda}+\bm{\delta})^{-1}
  \circ\mathcal{A}(\bm{\lambda})\circ\phi_0(x;\bm{\lambda})
  =B(0;\bm{\lambda})\varphi(x)^{-1}
  (1-e^{\partial}),\\
  \mathcal{B}(\bm{\lambda})
  &\eqdef\frac{1}{\sqrt{B(0;\bm{\lambda})}}\,
  \phi_0(x;\bm{\lambda})^{-1}\circ\mathcal{A}(\bm{\lambda})^{\dagger}
  \circ\phi_0(x;\bm{\lambda}+\bm{\delta})\n
  &=B(0;\bm{\lambda})^{-1}
  \bigl(B(x;\bm{\lambda})-D(x)e^{-\partial}\bigr)\varphi(x)\n
  &=\phi_0(x;\bm{\lambda})^{-2}\circ(1-e^{-\partial})\varphi(x)^{-1}\circ
  \phi_0(x;\bm{\lambda}+\bm{\delta})^2.
\end{align}
Starting from $\check{P}_0(x)=1$, $\check{P}_n(x;\bm{\lambda})$ can be
written as
\begin{align}
  \check{P}_n(x;\bm{\lambda})=P_n\bigl(\eta(x);\bm{\lambda}\bigr)
  &=\mathcal{B}(\bm{\lambda})
  \mathcal{B}(\bm{\lambda}+\bm{\delta})\cdots
  \mathcal{B}\bigl(\bm{\lambda}+(n-1)\bm{\delta}\bigr)\cdot 1\n
  &=\phi_0(x;\bm{\lambda})^{-2}\,
  \Bigl((1-e^{-\partial})\varphi(x)^{-1}\Bigr)^n\cdot
  \phi_0(x;\bm{\lambda}+n\bm{\delta})^2.
  \label{univrod}
\end{align}
This is the universal Rodrigues formula mentioned above.
Compare this with the individual Rodrigues formulas listed in \cite{kls}.
We note the following relations between the auxiliary function $\varphi(x)$
and the ground state eigenvectors, etc:
\begin{align}
  &\varphi(x)=\sqrt{\frac{B(0;\bm{\lambda})}{B(x;\bm{\lambda})}}
  \,\frac{\phi_0(x;\bm{\lambda}+\bm{\delta})}{\phi_0(x;\bm{\lambda})}
  =\frac{\eta(x+1)-\eta(x)}{\eta(1)},\\
  &\frac{B(x;\bm{\lambda}+\bm{\delta})}{B(x+1;\bm{\lambda})}
  =\kappa^{-1}\frac{\varphi(x+1)}{\varphi(x)}=\kappa^{-2}.
\end{align}

\subsection{Virtual systems}
\label{sec:virt}

\subsubsection{modified potential functions}
\label{sec:mpf}

Let us assume the existence of modified potential functions $B'(x)$, $D'(x)$
satisfying
\begin{alignat}{2}
  \alpha^2B'(x)D'(x+1)&=B(x)D(x+1),&\quad\ \alpha&>0,
  \label{B'D'=BD}\\
  \alpha\bigl(B'(x)+D'(x)\bigr)+\alpha'&=B(x)+D(x),&\alpha'&<0,
  \label{B'+D'=B+D}
\end{alignat}
where $\alpha$ and $\alpha'$ are constant.
Here $B'(x)$ and $D'(x)$ are rational functions of $x$ or $q^x$.
Note that these two equations are valid for any $x$($\in\mathbb{C}$).
These functions are required to satisfy the boundary conditions, too:
\begin{equation}
  B'(x)>0\ \ (x\in\mathbb{Z}_{\geq 0}),\quad
  D'(x)>0\ \ (x\in\mathbb{Z}_{\geq 1}),\ \ D'(0)=0.
  \label{I:B'>0}
\end{equation}
Then we obtain a linear relation between the original Hamiltonian
\eqref{Hdef2} and the virtual Hamiltonian $\mathcal{H}'$ \eqref{H'}
constructed by the modified potential functions $B'(x)$ and $D'(x)$,
\eqref{B'D'=BD}--\eqref{I:B'>0}:
\begin{align}
  \mathcal{H}&=\alpha\mathcal{H}'+\alpha',
  \label{HH'}\\
  \mathcal{H}'&\eqdef-\sqrt{B'(x)}\,e^{\partial}\sqrt{D'(x)}
  -\sqrt{D'(x)}\,e^{-\partial}\sqrt{B'(x)}+B'(x)+D'(x).
  \label{H'}
\end{align}
We define a positive function $\phi'_0(x)$ by
\begin{equation}
  \phi'_0(x)\eqdef\prod_{y=0}^{x-1}\sqrt{\frac{B'(y)}{D'(y+1)}}\,,
  \label{tphi0}
\end{equation}
which is annihilated by $\mathcal{A}'=\sqrt{B'(x)}-e^{\partial}\sqrt{D'(x)}$ :
\begin{equation}
  \mathcal{A}'\phi'_0(x)=0
  \ \Rightarrow \ \sqrt{B'(x)}\,\phi'_0(x)=\sqrt{D'(x+1)}\,\phi'_0(x).
\end{equation}
In other words, $\phi'_0(x)$ is a solution of the original Schr\"odinger
equation of $\mathcal{H}$ \eqref{schreq} with a negative energy $\alpha'$.
That is, $\phi'_0(x)$ does not belong to the $\ell^2$ Hilbert space spanned
by the eigenvectors of $\mathcal{H}$.
We introduce another positive function $\nu(x)$ by the ratio
$\phi_0(x)/\phi'_0(x)$:
\begin{equation}
  \nu(x)\eqdef\frac{\phi_0(x)}{\phi'_0(x)}
  =\prod_{y=0}^{x-1}\frac{B(y)}{\alpha B'(y)}
  =\prod_{y=0}^{x-1}\frac{\alpha D'(y+1)}{D(y+1)}.
  \label{phi0/tphi0}
\end{equation}
It can be analytically continued into a meromorphic function of $x$ or $q^x$
through the functional relations:
\begin{equation}
  \nu(x+1)=\frac{B(x)}{\alpha B'(x)}\nu(x),\quad
  \nu(x-1)=\frac{D(x)}{\alpha D'(x)}\nu(x).
  \label{nurel}
\end{equation}

As we will show shortly, these modified potential functions $B'(x)$, $D'(x)$
\eqref{B'D'=BD}-\eqref{B'+D'=B+D} are obtained from the original potential
functions $B(x)$, $D(x)$ \eqref{Hdef} by {\em twisting} the parameters,
which is an essential step for introducing the exceptional and multi-indexed
orthogonal polynomials \cite{os16,os25,os26,os27}.
Other solutions of the Schr\"odinger equation for the virtual Hamiltonian
$\mathcal{H}'$ can be obtained in a factorised form:
\begin{align}
  &\mathcal{H}'\tilde{\phi}_{\text{v}}(x)
  =\mathcal{E}'_{\text{v}}\tilde{\phi}_{\text{v}}(x)
  \ \ \Bigl(\Rightarrow\ \mathcal{H}\tilde{\phi}_{\text{v}}(x)
  =\tilde{\mathcal{E}}_{\text{v}}\tilde{\phi}_{\text{v}}(x),
  \ \ \tilde{\mathcal{E}}_{\text{v}}
  \eqdef\alpha\mathcal{E}'_{\text{v}}+\alpha'\Bigr),
  \label{Hphitv=}\\
  &\tilde{\phi}_{\text{v}}(x)\eqdef\phi'_0(x)\check{\xi}_{\text{v}}(x),
  \label{phitv}
\end{align}
in which $\check{\xi}_{\text{v}}(x)$ is a solution of
\begin{align}
  &\widetilde{\mathcal{H}}'\check{\xi}_{\text{v}}(x)
  =\mathcal{E}'_{\text{v}}\check{\xi}_{\text{v}}(x),\\
  &\widetilde{\mathcal{H}}'
  \eqdef\phi'_0(x)^{-1}\circ\mathcal{H}'\circ\phi'_0(x)
  =B'(x)(1-e^{\partial})+D'(x)(1-e^{-\partial}),
  \label{Htp}
\end{align}
namely,
\begin{equation}
  B'(x)\bigl(\check{\xi}_{\text{v}}(x)-\check{\xi}_{\text{v}}(x+1)\bigr)
  +D'(x)\bigl(\check{\xi}_{\text{v}}(x)-\check{\xi}_{\text{v}}(x-1)\bigr)
  =\mathcal{E}'_{\text{v}}\check{\xi}_{\text{v}}(x).
  \label{xieq}
\end{equation}
For a non-negative integer $\text{v}$, $\mathcal{E}'_{\text{v}}$ and
$\check{\xi}_{\text{v}}(x)$
can be obtained from the original polynomial solution
$\check{P}_{\text{v}}(x)=P_{\text{v}}(\eta(x))$ \eqref{tHcPn=} by twisting
the parameters \eqref{twist}.
Note that it also satisfies the universal normalisation condition
$\check{\xi}_{\text{v}}(0)=1$.
We call a solution
$\tilde{\phi}_{\text{v}}(x)=\phi'_0(x)\check{\xi}_{\text{v}}(x)$
{\em a virtual state vector}
if it has negative energy and positive on the entire semi-infinite lattice:
\begin{equation}
  \mathcal{E}'_{\text{v}}<0,\quad\check{\xi}_{\text{v}}(x)>0
  \ \ (x\in\mathbb{Z}_{\geq 0}).
  \label{xi>0}
\end{equation}
The positivity is necessary for the
well-definedness of the Darboux transformations in terms of
these virtual state vectors. See, for example, \eqref{BDd1def}.
Let us emphasise that the virtual state vector $\tilde{\phi}_{\text{v}}(x)$
is also a solution of the original Schr\"odinger equation \eqref{Hphitv=}
with a negative energy $\tilde{\mathcal{E}}_{\text{v}}<0$ and of infinite
norm $\norm{\tilde{\phi}_{\text{v}}}=\infty$.
The set $\mathcal{V}$ of non-negative integers $\text{v}$, satisfying the
above two conditions \eqref{xi>0} is called the index set of the virtual
state vectors.

\subsubsection{virtual state vectors}
\label{sec:vir}

In this subsection we present the explicit forms of the modified potential
functions $B'(x)$, $D'(x)$ \eqref{B'D'=BD}-\eqref{B'+D'=B+D} and the virtual
state vectors $\tilde{\phi}_{\text{v}}(x)=\phi'_0(x)\check{\xi}_{\text{v}}(x)$
\eqref{phitv} for M, l$q$J and l$q$L polynomials through twisting.

Let us define the twist operation $\mathfrak{t}$,
which is an {\em involution\/},
\begin{equation}
  \mathfrak{t}(\bm{\lambda})\eqdef\left\{
  \begin{array}{ll}
  (\lambda_1,\lambda_2^{-1})&:\text{M}\\
  (-\lambda_1,\lambda_2)&:\text{l$q$J}\\
  -\lambda_1&:\text{l$q$L}
  \end{array}\right.,\quad
  \text{namely}\ \left\{
  \begin{array}{ll}
  \mathfrak{t}(\bm{\lambda})=(\beta,c^{-1})&:\text{M}\\
  q^{\mathfrak{t}(\bm{\lambda})}=(a^{-1},b)&:\text{l$q$J}\\
  q^{\mathfrak{t}(\bm{\lambda})}=a^{-1}&:\text{l$q$L}
  \end{array}\right.,\quad
  \mathfrak{t}^2=\text{id},
  \label{twist}
\end{equation}
and the two functions $B'(x)$ and $D'(x)$,
\begin{equation}
  B'(x;\bm{\lambda})\eqdef B\bigl(x;\mathfrak{t}(\bm{\lambda})\bigr),\quad
  D'(x)\eqdef D(x),
  \label{B'D'}
\end{equation}
namely,
\begin{equation}
  B'(x;\bm{\lambda})=\left\{
  \begin{array}{ll}
  c^{-1}(x+\beta)&:\text{M}\\
  a^{-1}(q^{-x}-bq)&:\text{l$q$J}\\
  a^{-1}q^{-x}&:\text{l$q$L}
  \end{array}\right.,\quad
  D'(x)=\left\{
  \begin{array}{ll}
  x&:\text{M}\\
  q^{-x}-1&:\text{l$q$J,\,l$q$L}
  \end{array}\right..
\end{equation}
Then the conditions \eqref{B'D'=BD}--\eqref{I:B'>0} are satisfied with the
following $\alpha$ and $\alpha'$ :
\begin{equation}
  \alpha(\bm{\lambda})=\left\{
  \begin{array}{ll}
  c&:\text{M}\\
  a&:\text{l$q$J,\,l$q$L}
  \end{array}\right.,\quad
  \alpha'(\bm{\lambda})=\left\{
  \begin{array}{ll}
  -(1-c)\beta&:\text{M}\\
  -(1-a)(1-bq)&:\text{l$q$J}\\
  -(1-a)&:\text{l$q$L}
  \end{array}\right..
\end{equation}

The virtual Hamiltonian $\mathcal{H}'$, $\mathcal{E}'_{\text{v}}$ and
the virtual state vector $\tilde{\phi}_{\text{v}}(x)$ are given by
$\mathcal{H}(\mathfrak{t}(\bm{\lambda}))$,
$\mathcal{E}_{\text{v}}(\mathfrak{t}(\bm{\lambda}))$ and
$\phi_{\text{v}}(x;\mathfrak{t}(\bm{\lambda}))$. Namely,
\begin{align}
  &\mathcal{H}(\bm{\bm{\lambda}})=
  \alpha(\bm{\lambda})\mathcal{H}\bigl(\mathfrak{t}(\bm{\lambda})\bigr)
  +\alpha'(\bm{\lambda}),
  \label{H=aH'+a'}\\
  &\tilde{\phi}_{\text{v}}(x;\bm{\lambda})
  \eqdef\phi_{\text{v}}\bigl(x;\mathfrak{t}(\bm{\lambda})\bigr)
  =\tilde{\phi}_0(x;\bm{\lambda})
  \check{\xi}_{\text{v}}(x;\bm{\lambda})
  \ \ (\text{v}\in\mathcal{V}),
  \label{tphiv=}\\
  &\phi'_0(x;\bm{\lambda})=\tilde{\phi}_0(x;\bm{\lambda})
  \eqdef\phi_0\bigl(x;\mathfrak{t}(\bm{\lambda})\bigr),
  \ \check{\xi}_{\text{v}}(x;\bm{\lambda})\eqdef
  \check{P}_{\text{v}}\bigl(x;\mathfrak{t}(\bm{\lambda})\bigr),
  \ \check{\xi}_{\text{v}}(x;\bm{\lambda})\eqdef
  \xi_{\text{v}}\bigl(\eta(x);\bm{\lambda}\bigr),
  \label{xiv=}\\
  &\mathcal{H}(\bm{\lambda})\tilde{\phi}_{\text{v}}(x;\bm{\lambda})
  =\tilde{\mathcal{E}}_{\text{v}}(\bm{\lambda})
  \tilde{\phi}_{\text{v}}(x;\bm{\lambda}),\quad
  \mathcal{E}'_{\text{v}}(\bm{\lambda})
  =\mathcal{E}_{\text{v}}\bigl(\mathfrak{t}(\bm{\lambda})\bigr),\\
  &\tilde{\mathcal{E}}_{\text{v}}(\bm{\lambda})
  =\alpha(\bm{\lambda})\mathcal{E}'_{\text{v}}(\bm{\lambda})
  +\alpha'(\bm{\lambda})
  =\left\{
  \begin{array}{ll}
  -(1-c)(\text{v}+\beta)&:\text{M}\\
  -(1-aq^{-\text{v}})(1-bq^{\text{v}+1})&:\text{l$q$J}\\
  -(1-aq^{-\text{v}})&:\text{l$q$L}
  \end{array}\right.,
  \label{tEv=}\\
  &\nu(x;\bm{\lambda})=\left\{
  \begin{array}{ll}
  c^x&:\text{M}\\
  a^x&:\text{l$q$J,\,l$q$L}
  \end{array}\right..
\end{align}
Note that $\alpha'(\bm{\lambda})=\tilde{\mathcal{E}}_0(\bm{\lambda})<0$.
The negative virtual state energy condition \eqref{xi>0} is automatically
satisfied for M case but it restricts the maximal possible degree
$\text{v}_\text{max}$ for a given parameter $a$ for l$q$J and l$q$L cases as
\begin{equation}
  0<a<q^{\text{v}_\text{max}}.
  \label{rangelqJ}
\end{equation}

For the positivity of $\check{\xi}_{\text{v}}(x)$ \eqref{xi>0},
we write down $\check{\xi}_{\text{v}}(x)$ explicitly.
For M case, it is
\begin{align}
  \check{\xi}_{\text{v}}(x;\bm{\lambda})
  &={}_2F_1\Bigl(\genfrac{}{}{0pt}{}{-\text{v},\,-x}{\beta}\Bigm|1-c\Bigr)
  =\sum_{k=0}^{\min(\text{v},x)}
  \frac{(-\text{v},-x)_k}{(\beta)_k}\frac{(1-c)^k}{k!}\n
  &=\sum_{k=0}^{\min(\text{v},x)}\frac{(\text{v}-k+1,x-k+1)_k}{(\beta)_k}
  \frac{(1-c)^k}{k!},
  \label{Mvirt}
\end{align}
where we have used $(a)_k=(-1)^k(-a-k+1)_k$.
Since each $k$-th term of the sum in the last expression is positive,
the positivity of $\check{\xi}_{\text{v}}(x;\bm{\lambda})$ \eqref{xi>0} is
shown for all non-negative integer $\text{v}$.
For l$q$J case, by using the identity ((1.13.17) in \cite{kls})
\begin{equation*}
  {}_2\phi_1\Bigl(\genfrac{}{}{0pt}{}{q^{-n},\,b}{c}\Bigm|q;z\Bigr)
  =(bc^{-1}q^{-n}z;q)_n\,
  {}_3\phi_2\Bigl(\genfrac{}{}{0pt}{}{q^{-n},\,b^{-1}c,\,0}
  {c,\,b^{-1}cqz^{-1}}\Bigm|q;q\Bigr)\ \ (n\in\mathbb{Z}_{\geq 0}),
\end{equation*}
$\check{\xi}_{\text{v}}(x)$ is rewritten as
\begin{align}
  \check{\xi}_{\text{v}}(x;\bm{\lambda})
  &=\frac{(aq^{-\text{v}};q)_{\text{v}}}{(bq;q)_{\text{v}}}\,
  {}_2\phi_1\Bigl(\genfrac{}{}{0pt}{}{q^{-\text{v}},\,a^{-1}bq^{\text{v}+1}}
  {a^{-1}q}\Bigm|q;q^{x+1}\Bigr)\n
  &=\frac{(aq^{-\text{v}};q)_{\text{v}}}{(bq;q)_{\text{v}}}\,
  (bq^{x+1};q)_{\text{v}}\,
  {}_3\phi_2\Bigl(\genfrac{}{}{0pt}{}{q^{-\text{v}},\,b^{-1}q^{-\text{v}},\,0}
  {a^{-1}q,\,b^{-1}q^{-\text{v}-x}}\Bigm|q;q\Bigr)\n
  &=\frac{(aq^{-\text{v}},bq^{x+1};q)_{\text{v}}}{(bq;q)_{\text{v}}}
  \sum_{k=0}^{\text{v}}
  \frac{(q^{-\text{v}},b^{-1}q^{-\text{v}};q)_k}
  {(a^{-1}q,b^{-1}q^{-\text{v}-x};q)_k}\frac{q^k}{(q;q)_k}\n
  &=\frac{(aq^{-\text{v}},bq^{x+1};q)_{\text{v}}}{(bq;q)_{\text{v}}}
  \sum_{k=0}^{\text{v}}
  \frac{(q^{\text{v}-k+1},bq^{\text{v}-k+1};q)_k}
  {(aq^{-k},bq^{\text{v}-k+1+x},q;q)_k}(aq^{-\text{v}+x})^k,
  \label{lqJvirt}
\end{align}
where we have used $(a;q)_k=(-a)^kq^{\frac12k(k-1)}(a^{-1}q^{-k+1};q)_k$.
Since each $k$-th term of the sum and the overall factor in the last
expression is positive for $0<a<q^{\text{v}}$ \eqref{rangelqJ}, the
positivity of $\check{\xi}_{\text{v}}(x;\bm{\lambda})$ \eqref{xi>0} is shown.
By setting $b=0$, the same conclusion for l$q$L case is obtained.\\
{\bf Remark}\ For l$q$J, $\check{\xi}_{\text{v}}(x;\bm{\lambda})$ is
generically a degree $\text{v}$ polynomial in $\eta(x)=1-q^x$.
However, when $a=bq^{m+1}$, $\text{v}\leq m\in\mathbb{Z}_{\ge0}$,
$(a^{-1}bq^{\text{v}+1};q)_k=0$ for $k\ge m-\text{v}+1$.
That is, the highest degree is $m-\text{v}$.
To sum up, if $a=bq^{m+1}$ and $\text{v}\leq m<2\text{v}$,
$\check{\xi}_{\text{v}}(x;\bm{\lambda})$ is not a degree $\text{v}$ polynomial.
In later discussion we assume to avoid these special configurations of
the parameters.
\label{special}

The index set of the virtual state vectors is
\begin{equation}
  \mathcal{V}=\left\{
  \begin{array}{ll}
  \{1,2,\ldots\}&:\text{M}\\
  \{1,2,\ldots,\text{v}_{\text{max}}\}&:\text{l$q$J,\,l$q$L}
  \end{array}\right.,
  \label{vrange}
\end{equation}
where $\text{v}_{\text{max}}$ is the greatest integer satisfying
$0<a<q^{\text{v}}$.

The non square-summability of the virtual state vector
$\tilde{\phi}_{\text{v}}$ ($\norm{\tilde{\phi}_{\text{v}}}=\infty$) can be
directly verified through its large $x$ behaviour
\begin{align*}
  \tilde{\phi}_0(x;\bm{\lambda})^2\simeq\left\{
  \begin{array}{ll}
   \Gamma(\beta)^{-1}x^{\beta-1}c^{-x}&:\text{M}\\
   (bq;q)_{\infty}(q;q)_{\infty}^{-1}(a^{-1}q)^x&:\text{l$q$J}\\[2pt]
   (q;q)_{\infty}^{-1}(a^{-1}q)^x&:\text{l$q$L}
  \end{array}\right.\!\!\!,\quad
  \tilde{\xi}_{\text{v}}(x;\bm{\lambda})\simeq\left\{
  \begin{array}{ll}
   \text{const}\times x^{\text{v}}&:\text{M}\\
   \text{const}&:\text{l$q$J,\,l$q$L}
  \end{array}\right..
\end{align*}

Before closing this section, let us emphasise that the twisting is based
on the discrete symmetries of the difference Schr\"odinger equations
\eqref{schreq}, not of the equations governing the polynomials
\eqref{Hteq}-\eqref{tHcPn=}, as evidenced by the linear relation between
the two Hamiltonians \eqref{HH'}.

\section{Method of Virtual State Deletion}
\label{sec:mvsd}

We present the method of virtual states deletion for the discrete quantum
mechanics with real shifts (rdQM) in terms of the discrete QM analogue of
the multiple Darboux transformations \cite{os22,os26}.
The concept of the virtual states on a semi-infinite lattice is slightly
different from those on finite lattices as introduced in \cite{os26}.
On semi-infinite lattices, the solutions of the Schr\"odinger (eigenvalue)
equation \eqref{schreq} are eigenvectors with finite norm and the other
solutions with infinite norms, to which the virtual state vectors belong.
In contrast, the solutions of the Schr\"odinger (eigenvalue) equation on
finite lattices are eigenvectors only. This requires a subtle definition
of the virtual state vectors as introduced in \cite{os26}.

In the following discussion, the discrete analogue of the Wronskians,
to be called Casoratians, play an important role.
Here is a summary of useful properties of Casoratians.

\subsection{Casoratian}
\label{sec:Cas}

The Casorati determinant of a set of $n$ functions $\{f_j(x)\}$ is defined by
\begin{equation}
  \text{W}_{\text{C}}[f_1,\ldots,f_n](x)
  \eqdef\det\Bigl(f_k(x+j-1)\Bigr)_{1\leq j,k\leq n},
\end{equation}
(for $n=0$, we set $\text{W}_{\text{C}}[\cdot](x)=1$).
It satisfies identities ($n\geq 0$)
\begin{align}
  &\text{W}_{\text{C}}[gf_1,gf_2,\ldots,gf_n](x)
  =\prod_{k=0}^{n-1}g(x+k)\cdot\text{W}_{\text{C}}[f_1,f_2,\ldots,f_n](x),
  \label{Wformula1}\\
  &\text{W}_{\text{C}}\bigl[\text{W}_{\text{C}}[f_1,f_2,\ldots,f_n,g],
  \text{W}_{\text{C}}[f_1,f_2,\ldots,f_n,h]\,\bigr](x)\n
  &=\text{W}_{\text{C}}[f_1,f_2,\ldots,f_n](x+1)\,
  \text{W}_{\text{C}}[f_1,f_2,\ldots,f_n,g,h](x),
  \label{Wformula2}\\
  &\text{W}_{\text{C}}[F_1,F_2,\ldots,F_n](x)
  =(-1)^{\frac12n(n-1)}\prod_{k=0}^{n-2}
  \text{W}_{\text{C}}[f_1,f_2,\ldots,f_n](x+k),\n
  &\qquad\qquad\qquad\qquad\text{where }
  F_j(x)=\text{W}_{\text{C}}[f_1,\ldots,f_{j-1},f_{j+1},\ldots,f_n](x).
  \label{Wformula4}
\end{align}
The first and second identities will be used in the next subsection.
The third one is not explicitly used in this paper but it is used to express
the step backward shift operator in a determinant form \cite{rrmiop3}.

\subsection{Virtual states deletion}
\label{sec:vsd}

We delete virtual state vectors
$\tilde{\phi}_{d_1}(x),\tilde{\phi}_{d_2}(x),\ldots$
($d_j$\ :\ mutually distinct)\footnote{
Although this notation $d_j$ conflicts with the notation of the normalisation
constant $d_n$ in \eqref{ortho}, we think this does not cause any confusion
because the latter appears as $\frac{1}{d_n^2}\,\delta_{nm}$.
}
by the multiple application of Darboux transformations.

\subsubsection{one virtual state vector deletion}
\label{sec:one}

First we rewrite the original Hamiltonian $\mathcal{H}$ \eqref{Hdef} in a
factorised form in which the first virtual state vector
$\tilde{\phi}_{d_1}(x)$ ($d_1\in\mathcal{V}$) is annihilated by
$\hat{\mathcal{A}}_{d_1}$,
$\hat{\mathcal{A}}_{d_1}\tilde{\phi}_{d_1}(x)=0$:
\begin{align*}
  &\mathcal{H}=\hat{\mathcal{A}}_{d_1}^{\dagger}\hat{\mathcal{A}}_{d_1}
  +\tilde{\mathcal{E}}_{d_1},\\
  &\hat{\mathcal{A}}_{d_1}\eqdef
  \sqrt{\hat{B}_{d_1}(x)}-e^{\partial}\sqrt{\hat{D}_{d_1}(x)},\quad
  \hat{\mathcal{A}}_{d_1}^{\dagger}
  =\sqrt{\hat{B}_{d_1}(x)}-\sqrt{\hat{D}_{d_1}(x)}\,e^{-\partial}.
\end{align*}
Here the potential functions
$\hat{B}_{d_1}(x)$ and $\hat{D}_{d_1}(x)$ are determined by $B'(x)$, $D'(x)$
and the virtual state polynomial $\check{\xi}_{d_1}(x)$:
\begin{equation}
  \hat{B}_{d_1}(x)\eqdef\alpha B'(x)
  \frac{\check{\xi}_{d_1}(x+1)}{\check{\xi}_{d_1}(x)},\quad
  \hat{D}_{d_1}(x)\eqdef\alpha D'(x)
  \frac{\check{\xi}_{d_1}(x-1)}{\check{\xi}_{d_1}(x)}.
  \label{BDd1def}
\end{equation}
We have $\hat{B}_{d_1}(x)>0$ ($x\in\mathbb{Z}_{\geq 0}$),
$\hat{D}_{d_1}(x)>0$ ($x\in\mathbb{Z}_{\geq 1}$), $\hat{D}_{d_1}(0)=0$ and
\begin{align*}
  \hat{B}_{d_1}(x)\hat{D}_{d_1}(x+1)&=B(x)D(x+1),\\
  \hat{B}_{d_1}(x)+\hat{D}_{d_1}(x)+\tilde{\mathcal{E}}_{d_1}&=B(x)+D(x),
\end{align*}
where use is made of \eqref{xieq} in the second equation.

Next let us define a new Hamiltonian $\mathcal{H}_{d_1}$ by changing
the order of the two matrices $\hat{\mathcal{A}}_{d_1}^{\dagger}$ and
$\hat{\mathcal{A}}_{d_1}$ together with the sets of new eigenvectors
$\phi_{d_1n}(x)$ and new virtual state vectors
$\tilde{\phi}_{d_1\text{v}}(x)$:
\begin{align}
  \mathcal{H}_{d_1}&\eqdef
  \hat{\mathcal{A}}_{d_1}\hat{\mathcal{A}}_{d_1}^{\dagger}
  +\tilde{\mathcal{E}}_{d_1},\quad
  \mathcal{H}_{d_1}=(\mathcal{H}_{d_1\,x,y})
  \ \ (x,y\in\mathbb{Z}_{\geq 0}),\\
  \phi_{d_1n}(x)&\eqdef\hat{\mathcal{A}}_{d_1}\phi_n(x)
  \ \ (n\in\mathbb{Z}_{\geq 0}),
  \label{phid1n}\\
  \tilde{\phi}_{d_1\text{v}}(x)&\eqdef
  \hat{\mathcal{A}}_{d_1}\tilde{\phi}_{\text{v}}(x)
  \ \ (\text{v}\in\mathcal{V}\backslash\{d_1\}).
  \label{tphid1v}
\end{align}
It is easy to verify that $\phi_{d_1n}(x)$ is an eigenvector and that
$\tilde{\phi}_{d_1\text{v}}(x)$ is a virtual state vector
\begin{equation*}
  \mathcal{H}_{d_1}\phi_{d_1n}(x)
  =\mathcal{E}_n\phi_{d_1n}(x)
  \ \ (n\in\mathbb{Z}_{\geq 0}),\quad
  \mathcal{H}_{d_1}\tilde{\phi}_{d_1\text{v}}(x)
  =\tilde{\mathcal{E}}_{\text{v}}\tilde{\phi}_{d_1\text{v}}(x)
  \ \ (\text{v}\in\mathcal{V}\backslash\{d_1\}).
\end{equation*}
For example,
\begin{align*}
  &\mathcal{H}_{d_1}\phi_{d_1n}
  =(\hat{\mathcal{A}}_{d_1}\hat{\mathcal{A}}_{d_1}^{\dagger}
  +\tilde{\mathcal{E}}_{d_1})\hat{\mathcal{A}}_{d_1}\phi_n
  =\hat{\mathcal{A}}_{d_1}(\hat{\mathcal{A}}_{d_1}^{\dagger}
  \hat{\mathcal{A}}_{d_1}+\tilde{\mathcal{E}}_{d_1})\phi_n\\
  &=\hat{\mathcal{A}}_{d_1}\mathcal{H}\phi_n
  =\hat{\mathcal{A}}_{d_1}\mathcal{E}_n\phi_n
  =\mathcal{E}_n\hat{\mathcal{A}}_{d_1}\phi_n
  =\mathcal{E}_n\phi_{d_1n}.
\end{align*}
The two Hamiltonians $\mathcal{H}$ and $\mathcal{H}_{d_1}$ are {\em exactly
iso-spectral}. If the original system is exactly solvable, this new system
is also exactly solvable.
The orthogonality relation for the new eigenvectors is
\begin{align*}
  (\phi_{d_1n},\phi_{d_1m})
  &=(\hat{\mathcal{A}}_{d_1}\phi_n,\hat{\mathcal{A}}_{d_1}\phi_m)
  =(\hat{\mathcal{A}}_{d_1}^{\dagger}\hat{\mathcal{A}}_{d_1}\phi_n,\phi_m)
  =\bigl((\mathcal{H}-\tilde{\mathcal{E}}_{d_1})\phi_n,\phi_m\bigr)\n
  &=(\mathcal{E}_n-\tilde{\mathcal{E}}_{d_1})(\phi_n,\phi_m)
  =(\mathcal{E}_n-\tilde{\mathcal{E}}_{d_1})\frac{\delta_{nm}}{d_n^2}
  \ \ (n,m\in\mathbb{Z}_{\geq 0}).
\end{align*}
The second equality requires the condition
$\lim\limits_{N\to\infty}\sqrt{\hat{D}_{d_1}(N+1)}\,
\phi_{d_1\,n}(N)\phi_m(N+1)=0$,
which is easily verified for the three systems
(M, l$q$J and l$q$L).
This shows clearly that the {\em negative} virtual state energy
($\tilde{\mathcal{E}}_{\text{v}}<0$) is necessary for the positivity of
the inner products.

The new eigenvector $\phi_{d_1n}(x)$ \eqref{phid1n} and the virtual state
vector $\tilde{\phi}_{d_1\text{v}}(x)$ \eqref{tphid1v} are expressed neatly
in terms of the Casoratians:
\begin{equation}
  \phi_{d_1n}(x)
  =\frac{-\sqrt{\alpha B'(x)}\,\phi'_0(x)
  \text{W}_{\text{C}}[\check{\xi}_{d_1},\nu\check{P}_n](x)}
  {\sqrt{\check{\xi}_{d_1}(x)\check{\xi}_{d_1}(x+1)}},
  \ \ \tilde{\phi}_{d_1\text{v}}(x)
  =\frac{-\sqrt{\alpha B'(x)}\,\phi'_0(x)
  \text{W}_{\text{C}}[\check{\xi}_{d_1},\check{\xi}_{\text{v}}](x)}
  {\sqrt{\check{\xi}_{d_1}(x)\check{\xi}_{d_1}(x+1)}}.
  \label{phid1nW}
\end{equation}
The positivity of the virtual state vectors $\tilde{\phi}_{\text{v}}(x)$
is inherited by the new virtual state vectors $\tilde{\phi}_{d_1\text{v}}(x)$
\eqref{tphid1v}.
The Casoratian
$\text{W}_{\text{C}}[\check{\xi}_{d_1},\check{\xi}_{\text{v}}](x)$
has definite sign for $x\in\mathbb{Z}_{\geq 0}$, namely all positive
or all negative.
By using \eqref{xieq} we have
\begin{equation*}
  \alpha B'(x)\text{W}_{\text{C}}[\check{\xi}_{d_1},\check{\xi}_{\text{v}}](x)
  =\alpha D'(x)\text{W}_{\text{C}}[\check{\xi}_{d_1},\check{\xi}_{\text{v}}]
  (x-1)
  +(\tilde{\mathcal{E}}_{d_1}-\tilde{\mathcal{E}}_{\text{v}})
  \check{\xi}_{d_1}(x)\check{\xi}_{\text{v}}(x).
\end{equation*}
At $x=0$, $\text{W}_{\text{C}}[\check{\xi}_{d_1},\check{\xi}_{\text{v}}](0)$
has the same sign as
$(\tilde{\mathcal{E}}_{d_1}-\tilde{\mathcal{E}}_{\text{v}})$,
as $D'(0)=0$, $\alpha B'(0)>0$ and
$\check{\xi}_{d_1}(0)\check{\xi}_{\text{v}}(0)=1$.
By setting $x=1,2,\ldots$ in turn, we obtain \label{samesign}
\begin{equation*}
  \text{sgn}(\tilde{\mathcal{E}}_{d_1}-\tilde{\mathcal{E}}_{\text{v}})
  \text{W}_{\text{C}}[\check{\xi}_{d_1},\check{\xi}_{\text{v}}](x)>0
  \ \ (x\in\mathbb{Z}_{\geq 0}).
\end{equation*}
The new ground state eigenvector $\phi_{d_10}(x)$ is of
definite sign as the original one $\phi_0(x)$ \eqref{phi0}.
We show that the Casoratian $\text{W}_{\text{C}}[\check{\xi}_{d_1},\nu](x)$
has definite sign for $x\in\mathbb{Z}_{\geq 0}$.
By writing down the equation $\mathcal{H}\phi_n(x)=\mathcal{E}_n\phi_n(x)$
with $\mathcal{H}=\hat{\mathcal{A}}_{d_1}^{\dagger}\hat{\mathcal{A}}_{d_1}
+\tilde{\mathcal{E}}_{d_1}$, we obtain
\begin{equation}
  \alpha B'(x)\bigl(\nu\check{P}_n\bigr)(x+1)
  +\alpha D'(x)\bigl(\nu\check{P}_n\bigr)(x-1)
  =\bigl(\alpha B'(x)+\alpha D'(x)+\alpha'-\mathcal{E}_n\bigr)
  \bigl(\nu\check{P}_n\bigr)(x),
  \label{eqnuPn}
\end{equation}
where $\bigl(fg\bigr)(x)\eqdef f(x)g(x)$.
By using this, we can show
\begin{equation*}
  \alpha B'(x)\text{W}_{\text{C}}[\check{\xi}_{d_1},\nu\check{P}_n](x)
  =\alpha D'(x)\text{W}_{\text{C}}[\check{\xi}_{d_1},\nu\check{P}_n](x-1)
  +(\tilde{\mathcal{E}}_{d_1}-\mathcal{E}_n)
  \check{\xi}_{d_1}(x)\nu(x)\check{P}_n(x).
\end{equation*}
By setting $n=0$ and $x=0,1,2,\ldots$ in turn, we obtain
\begin{equation*}
  -\text{W}_{\text{C}}[\check{\xi}_{d_1},\nu](x)>0
  \ \ (x\in\mathbb{Z}_{\geq 0}).
\end{equation*}

Let us rewrite the deformed Hamiltonian $\mathcal{H}_{d_1}$ in the
standard form based on the ground state eigenvector $\phi_{d_10}(x)$.
Let us introduce the potential functions $B_{d_1}(x)$ and $D_{d_1}(x)$ by
\begin{align}
  B_{d_1}(x)&\eqdef\alpha B'(x+1)
  \frac{\check{\xi}_{d_1}(x)}{\check{\xi}_{d_1}(x+1)}
  \frac{\text{W}_{\text{C}}[\check{\xi}_{d_1},\nu](x+1)}
  {\text{W}_{\text{C}}[\check{\xi}_{d_1},\nu](x)},\\
  D_{d_1}(x)&\eqdef\alpha D'(x)
  \frac{\check{\xi}_{d_1}(x+1)}{\check{\xi}_{d_1}(x)}
  \frac{\text{W}_{\text{C}}[\check{\xi}_{d_1},\nu](x-1)}
  {\text{W}_{\text{C}}[\check{\xi}_{d_1},\nu](x)}.
\end{align}
The positivity of $B_{d_1}(x)$ and $D_{d_1}(x)$ is obvious and
the boundary condition $D_{d_1}(0)=0$ is satisfied.
They satisfy the relations
\begin{align*}
  B_{d_1}(x)D_{d_1}(x+1)
  &=\hat{B}_{d_1}(x+1)\hat{D}_{d_1}(x+1),\\
  B_{d_1}(x)+D_{d_1}(x)
  &=\hat{B}_{d_1}(x)+\hat{D}_{d_1}(x+1)+\tilde{\mathcal{E}}_{d_1}.
\end{align*}
The standard form Hamiltonian $\mathcal{H}_{d_1}$ is obtained:
\begin{align}
  &\mathcal{H}_{d_1}=\mathcal{A}_{d_1}^{\dagger}\mathcal{A}_{d_1},\\
  &\mathcal{A}_{d_1}\eqdef
  \sqrt{B_{d_1}(x)}-e^{\partial}\sqrt{D_{d_1}(x)},\quad
  \mathcal{A}_{d_1}^{\dagger}
  =\sqrt{B_{d_1}(x)}-\sqrt{D_{d_1}(x)}\,e^{-\partial},
\end{align}
in which $\mathcal{A}_{d_1}$ annihilates the ground state eigenvector,
$\mathcal{A}_{d_1}\phi_{d_10}(x)=0$.

It is worthwhile to emphasise that the new Hamiltonian $\mathcal{H}_{d_1}$
is constructed by using the virtual state vector $\tilde{\phi}_{d_1}(x)$,
which is now deleted from the new index set of available virtual state
vectors $\mathcal{V}\backslash\{d_1\}$.

\subsubsection{multi virtual state vectors deletion}
\label{sec:I:multi}

We repeat the above procedure and obtain the modified systems.
The number of deleted virtual state vectors should be less than or equal
$|\mathcal{V}|$.

For simplicity in notation, we introduce the following quantities ($s\geq 0$):
\begin{align}
  w_s(x)&\eqdef\text{W}_{\text{C}}[\check{\xi}_{d_1},\ldots,
  \check{\xi}_{d_s}](x),
  \label{ws}\\
  w'_{s,\text{v}}(x)&\eqdef\text{W}_{\text{C}}[\check{\xi}_{d_1},\ldots,
  \check{\xi}_{d_s},\check{\xi}_{\text{v}}](x),
  \label{w'sv}\\
  w''_{s,n}(x)&\eqdef\text{W}_{\text{C}}[\check{\xi}_{d_1},\ldots,
  \check{\xi}_{d_s},\nu\check{P}_n](x).
  \label{w''sn}
\end{align}
Note that $w'_{s,d_{s+1}}(x)=w_{s+1}(x)$.

Let us assume that we have already deleted $s$ virtual state vectors
($s\geq 1$), which are labeled by $\{d_1,\ldots,d_s\}$
($d_j\in\mathcal{V}$ : mutually distinct).
Namely we have
\begin{align}
  &\mathcal{H}_{d_1\ldots d_s}\eqdef
  \hat{\mathcal{A}}_{d_1\ldots d_s}\hat{\mathcal{A}}_{d_1\ldots d_s}^{\dagger}
  +\tilde{\mathcal{E}}_{d_s},\quad
  \mathcal{H}_{d_1\ldots d_s}=(\mathcal{H}_{d_1\ldots d_s\,x,y})\quad
  (x,y\in\mathbb{Z}_{\geq 0}),
  \label{Hd1..ds}\\
  &\hat{\mathcal{A}}_{d_1\ldots d_s}\eqdef
  \sqrt{\hat{B}_{d_1\dots d_s}(x)}
  -e^{\partial}\sqrt{\hat{D}_{d_1\ldots d_s}(x)},
  \quad\hat{\mathcal{A}}_{d_1\ldots d_s}^{\dagger}
  =\sqrt{\hat{B}_{d_1\ldots d_s}(x)}
  -\sqrt{\hat{D}_{d_1\ldots d_s}(x)}\,e^{-\partial},\!\!\\
  &\hat{B}_{d_1\ldots d_s}(x)\eqdef\alpha B'(x+s-1)
  \frac{w_{s-1}(x)}{w_{s-1}(x+1)}\frac{w_s(x+1)}{w_s(x)},
  \label{Bdsform}\\
  &\hat{D}_{d_1\ldots d_s}(x)\eqdef\alpha D'(x)
  \frac{w_{s-1}(x+1)}{w_{s-1}(x)}\frac{w_s(x-1)}{w_s(x)},
  \label{Ddsform}\\
  &\phi_{d_1\ldots d_s\,n}(x)\eqdef
  \hat{\mathcal{A}}_{d_1\ldots d_s}\phi_{d_1\ldots d_{s-1}\,n}(x)
  \ \ (n\in\mathbb{Z}_{\geq 0}),\\
  &\tilde{\phi}_{d_1\ldots d_s\,\text{v}}(x)\eqdef
  \hat{\mathcal{A}}_{d_1\ldots d_s}
  \tilde{\phi}_{d_1\ldots d_{s-1}\,\text{v}}(x)
  \ \ (\text{v}\in\mathcal{V}\backslash\{d_1,\ldots,d_s\}),\\
  &\mathcal{H}_{d_1\ldots d_s}\phi_{d_1\ldots d_s\,n}(x)
  =\mathcal{E}_n\phi_{d_1\ldots d_s\,n}(x)
  \ \ (n\in\mathbb{Z}_{\geq 0}),
  \label{Hd1..dsphid1..ds=}\\
  &\mathcal{H}_{d_1\ldots d_s}\tilde{\phi}_{d_1\ldots d_s\,\text{v}}(x)
  =\tilde{\mathcal{E}}_\text{v}\tilde{\phi}_{d_1\ldots d_s\,\text{v}}(x)
  \ \ (\text{v}\in\mathcal{V}\backslash\{d_1,\ldots,d_s\}),
  \label{Hd1..dstphid1..ds=}\\
  &(\phi_{d_1\ldots d_s\,n},\phi_{d_1\ldots d_s\,m})
  =\prod_{j=1}^s(\mathcal{E}_n-\tilde{\mathcal{E}}_{d_j})\cdot
  \frac{\delta_{nm}}{d_n^2}
  \ \ (n,m\in\mathbb{Z}_{\geq 0}).
  \label{(phid1..dsm,phid1..dsn)}
\end{align}
The eigenvectors and the virtual state vectors have Casoratian expressions:
\begin{align}
  \phi_{d_1\ldots d_s\,n}(x)
  &=\mathcal{S}_{d_1\ldots d_s}
  \frac{\sqrt{\prod_{j=1}^s\alpha B'(x+j-1)}\,\phi'_0(x)}
  {\sqrt{w_s(x)w_s(x+1)}}\,w''_{s,n}(x),
  \label{phid1..dsn}\\[4pt]
  \tilde{\phi}_{d_1\ldots d_s\,\text{v}}(x)
  &=\mathcal{S}_{d_1\ldots d_s}
  \frac{\sqrt{\prod_{j=1}^s\alpha B'(x+j-1)}\,\phi'_0(x)}
  {\sqrt{w_s(x)w_s(x+1)}}\,w'_{s,\text{v}}(x),
  \label{phitd1..dsv}
\end{align}
where the sign factor $\mathcal{S}_{d_1\ldots d_s}$ is
\begin{equation}
  \mathcal{S}_{d_1\ldots d_s}=(-1)^s
  \!\!\prod_{1\leq i<j\leq s}\!\!\text{sgn}(\tilde{\mathcal{E}}_{d_i}
  -\tilde{\mathcal{E}}_{d_j}).
  \label{S}
\end{equation}
The Casoratians of multiple virtual state vectors $w_s(x)$,
$w'_{s,\text{v}}(x)$ and with the ground state eigenvector
$w''_{s,0}(x)$ \eqref{ws}--\eqref{w''sn} are of definite sign
\begin{equation}
  \left.
  \begin{array}{r}
  \prod\limits_{1\leq i<j\leq s}\text{sgn}(\tilde{\mathcal{E}}_{d_i}
  -\tilde{\mathcal{E}}_{d_j})\cdot w_s(x)>0\\[4pt]
  \prod\limits_{1\leq i<j\leq s}\text{sgn}(\tilde{\mathcal{E}}_{d_i}
  -\tilde{\mathcal{E}}_{d_j})\cdot
  \prod\limits_{i=1}^s\text{sgn}(\tilde{\mathcal{E}}_{d_i}
  -\tilde{\mathcal{E}}_{\text{v}})\cdot w'_{s,\text{v}}(x)>0\\[4pt]
  \prod\limits_{1\leq i<j\leq s}\text{sgn}(\tilde{\mathcal{E}}_{d_i}
  -\tilde{\mathcal{E}}_{d_j})\cdot(-1)^sw''_{s,0}(x)>0
  \end{array}
  \right\}\ \ (x\in\mathbb{Z}_{\geq 0}).
  \label{wts>0}
\end{equation}
So we have $\hat{B}_{d_1\ldots d_s}(x)>0$ ($x\in\mathbb{Z}_{\geq 0}$),
$\hat{D}_{d_1\ldots d_s}(x)>0$ ($x\in\mathbb{Z}_{\geq 1}$)
and $\hat{D}_{d_1\ldots d_s}(0)=0$.
For \eqref{(phid1..dsm,phid1..dsn)}, the asymptotic condition
\begin{equation}
  \lim_{N\to\infty}\sqrt{\hat{D}_{d_1\ldots d_s}(N+1)}\,
  \phi_{d_1\ldots d_s\,n}(N)\phi_{d_1\ldots d_{s-1}\,m}(N+1)=0,
  \label{condAd}
\end{equation}
is necessary.
This can be easily verified for M, l$q$J and l$q$L systems.

By writing down \eqref{Hd1..dstphid1..ds=} and \eqref{Hd1..dsphid1..ds=}
explicitly, we obtain the following identities
\begin{align}
  &\quad\hat{B}_{d_1\ldots d_s}(x)+\hat{D}_{d_1\ldots d_s}(x+1)
  +\tilde{\mathcal{E}}_{d_s}-\tilde{\mathcal{E}}_{\text{v}}\n
  &=\alpha B'(x+s)\frac{w_s(x)}{w_s(x+1)}
  \frac{w'_{s,\text{v}}(x+1)}{w'_{s,\text{v}}(x)}
  +\alpha D'(x)\frac{w_s(x+1)}{w_s(x)}
  \frac{w'_{s,\text{v}}(x-1)}{w'_{s,\text{v}}(x)},
  \label{I:id1}\\
  &\quad\hat{B}_{d_1\ldots d_s}(x)+\hat{D}_{d_1\ldots d_s}(x+1)
  +\tilde{\mathcal{E}}_{d_s}-\mathcal{E}_n\n
  &=\alpha B'(x+s)\frac{w_s(x)}{w_s(x+1)}
  \frac{w''_{s,n}(x+1)}{w''_{s,n}(x)}
  +\alpha D'(x)\frac{w_s(x+1)}{w_s(x)}
  \frac{w''_{s,n}(x-1)}{w''_{s,n}(x)},
  \label{I:id2}
\end{align}
which are valid for any $x$($\in\mathbb{C}$).
Note that \eqref{I:id1}-\eqref{I:id2} are valid for $s=0$ by replacing
$\hat{B}_{d_1\ldots d_s}(x)+\hat{D}_{d_1\ldots d_s}(x+1)
+\tilde{\mathcal{E}}_{d_s}$ with $\alpha B'(x)+\alpha D'(x)+\alpha'$.
The identity \eqref{Wformula2} gives
\begin{equation}
  w_s(x+1)w''_{s+1,n}(x)=\text{W}_{\text{C}}[w_{s+1},w''_{s,n}](x),\quad
  w_s(x+1)w'_{s+1,\text{v}}(x)
  =\text{W}_{\text{C}}[w_{s+1},w'_{s,\text{v}}](x).
  \label{I:ww=Www}
\end{equation}
By using \eqref{I:id1}--\eqref{I:ww=Www},
we obtain the following identities ($s\geq 0$)
\begin{align}
  &\quad\alpha B'(x+s)w_s(x)w'_{s+1,\text{v}}(x)\n
  &=\alpha D'(x)w_s(x+1)w'_{s+1,\text{v}}(x-1)
  +(\tilde{\mathcal{E}}_{d_{s+1}}-\tilde{\mathcal{E}}_{\text{v}})
  w_{s+1}(x)w'_{s,\text{v}}(x),
  \label{I:id3}\\
  &\quad\alpha B'(x+s)w_s(x)w''_{s+1,n}(x)\n
  &=\alpha D'(x)w_s(x+1)w'_{s+1,n}(x-1)
  +(\tilde{\mathcal{E}}_{d_{s+1}}-\mathcal{E}_n)
  w_{s+1}(x)w''_{s,n}(x),
  \label{I:id4}
\end{align}
which are valid for any $x$($\in\mathbb{C}$).

The next step begins with rewriting the Hamiltonian
$\mathcal{H}_{d_1\ldots d_s}$ by choosing the next virtual state to be
deleted $d_{s+1}\in\mathcal{V}\backslash\{d_1,\ldots,d_s\}$.
The new potential functions
$\hat{B}_{d_1\ldots d_{s+1}}(x)$ and $\hat{D}_{d_1\ldots d_{s+1}}(x)$
are defined as in \eqref{Bdsform}--\eqref{Ddsform} by $s\to s+1$, 
and they contain the Casoratians \eqref{ws}--\eqref{w'sv} as ratios
(recall $w'_{s,d_{s+1}}(x)=w_{s+1}(x)$).
These Casoratians are of the same sign as shown in \eqref{wts>0},
the new potential functions are positive $\hat{B}_{d_1\ldots d_{s+1}}(x)>0$
($x\in\mathbb{Z}_{\geq 0}$),
$\hat{D}_{d_1\ldots d_{s+1}}(x)>0$ ($x\in\mathbb{Z}_{\geq 1}$) and
$\hat{D}_{d_1\ldots d_{s+1}}(0)=0$.
By using \eqref{I:id1}, we can show the relations
\begin{align*}
  \hat{B}_{d_1\ldots d_{s+1}}(x)\hat{D}_{d_1\ldots d_{s+1}}(x+1)
  &=\hat{B}_{d_1\ldots d_s}(x+1)\hat{D}_{d_1\ldots d_s}(x+1),\\
  \hat{B}_{d_1\ldots d_{s+1}}(x)+\hat{D}_{d_1\ldots d_{s+1}}(x)
  +\tilde{\mathcal{E}}_{d_{s+1}}
  &=\hat{B}_{d_1\ldots d_s}(x)+\hat{D}_{d_1\ldots d_s}(x+1)
  +\tilde{\mathcal{E}}_{d_s}.
\end{align*}
Therefore the Hamiltonian $\mathcal{H}_{d_1\ldots d_s}$ is rewritten as:
\begin{align*}
  &\mathcal{H}_{d_1\ldots d_s}
  =\hat{\mathcal{A}}_{d_1\ldots d_{s+1}}^{\dagger}
  \hat{\mathcal{A}}_{d_1\ldots d_{s+1}}
  +\tilde{\mathcal{E}}_{d_{s+1}},\\
  &\hat{\mathcal{A}}_{d_1\ldots d_{s+1}}\!\eqdef\!
  \sqrt{\hat{B}_{d_1\ldots d_{s+1}}(x)}
  -e^{\partial}\sqrt{\hat{D}_{d_1\ldots d_{s+1}}(x)},
  \ \hat{\mathcal{A}}_{d_1\ldots d_{s+1}}^{\dagger}\!
  =\!\sqrt{\hat{B}_{d_1\ldots d_{s+1}}(x)}
  -\sqrt{\hat{D}_{d_1\ldots d_{s+1}}(x)}\,e^{-\partial}.
\end{align*}

Now let us define a new Hamiltonian $\mathcal{H}_{d_1\ldots d_{s+1}}$
by changing the orders of
$\hat{\mathcal{A}}_{d_1\ldots d_{s+1}}^{\dagger}$ and
$\hat{\mathcal{A}}_{d_1\ldots d_{s+1}}$ together with
the eigenvectors $\phi_{d_1\ldots d_{s+1}\,n}(x)$ and
the virtual state vectors $\tilde{\phi}_{d_1\ldots d_{s+1}\,\text{v}}(x)$:
\begin{align*}
  &\mathcal{H}_{d_1\ldots d_{s+1}}\eqdef
  \hat{\mathcal{A}}_{d_1\ldots d_{s+1}}
  \hat{\mathcal{A}}_{d_1\ldots d_{s+1}}^{\dagger}
  +\tilde{\mathcal{E}}_{d_{s+1}},\quad
  \mathcal{H}_{d_1\ldots d_{s+1}}=(\mathcal{H}_{d_1\ldots d_{s+1}\,x,y})
  \ \ (x,y\in\mathbb{Z}_{\geq 0}),\\
  &\phi_{d_1\ldots d_{s+1}\,n}(x)\eqdef
  \hat{\mathcal{A}}_{d_1\ldots d_{s+1}}\phi_{d_1\ldots d_s\,n}(x)
  \ \ (n\in\mathbb{Z}_{\geq 0}),\\
  &\tilde{\phi}_{d_1\ldots d_{s+1}\,\text{v}}(x)\eqdef
  \hat{\mathcal{A}}_{d_1\ldots d_{s+1}}
  \tilde{\phi}_{d_1\ldots d_s\,\text{v}}(x)
  \ \ (\text{v}\in\mathcal{V}\backslash\{d_1,\ldots,d_{s+1}\}).
\end{align*}
It is easy to show that they satisfy
\eqref{Hd1..dsphid1..ds=}--\eqref{(phid1..dsm,phid1..dsn)} with $s\to s+1$.
By using \eqref{I:ww=Www}, we can show that the functions
$\phi_{d_1\ldots d_{s+1}\,n}(x)$ and
$\tilde{\phi}_{d_1\ldots d_{s+1}\,\text{v}}(x)$ are expressed as
\eqref{phid1..dsn}--\eqref{phitd1..dsv} with $s\to s+1$ and
$\mathcal{S}_{d_1\ldots d_{s+1}}$ satisfies
\begin{equation}
  \mathcal{S}_{d_1\ldots d_{s+1}}
  =-\mathcal{S}_{d_1\ldots d_s}
  \prod_{i=1}^s\text{sgn}(\tilde{\mathcal{E}}_{d_i}
  -\tilde{\mathcal{E}}_{d_{s+1}}),
  \label{Sd1..ds+1=}
\end{equation}
which is consistent with \eqref{S} and the initial value
$\mathcal{S}_{d_1}=-1$, see \eqref{phid1nW}.

The identities \eqref{I:id3}--\eqref{I:id4} will be used to show that the 
Casoratians $w'_{s+1,\text{v}}(x)$ and $w''_{s+1,0}(x)$ do not change sign for
$x\in\mathbb{Z}_{\geq 0}$, in the same manner as below \eqref{phid1nW}.
These establish the $s+1$ case.

At the end of this subsection we present this deformed Hamiltonian
$\mathcal{H}_{d_1\ldots d_s}$ in the standard form, in which the
$\mathcal{A}$ operator annihilates the ground state eigenvector:
\begin{align}
  &\mathcal{H}_{d_1\ldots d_s}
  =\mathcal{A}_{d_1\ldots d_s}^{\dagger}\mathcal{A}_{d_1\ldots d_s},
  \label{Hd1..ds=AdA}\\
  &\mathcal{A}_{d_1\ldots d_s}\eqdef
  \sqrt{B_{d_1\ldots d_s}(x)}
  -e^{\partial}\sqrt{D_{d_1\ldots d_s}(x)},
  \ \ \mathcal{A}_{d_1\ldots d_s}^{\dagger}
  =\sqrt{B_{d_1\ldots d_s}(x)}
  -\sqrt{D_{d_1\ldots d_s}(x)}\,e^{-\partial},
  \label{Ad1..ds=}
\end{align}
which satisfies $\mathcal{A}_{d_1\ldots d_s}\phi_{d_1\ldots d_s\,0}(x)=0$.
The potential functions $B_{d_1\ldots d_s}(x)$ and
$D_{d_1\ldots d_s}(x)$ are:
\begin{align}
  B_{d_1\ldots d_s}(x)&\eqdef\alpha B'(x+s)
  \frac{w_s(x)}{w_s(x+1)}\frac{w''_{s,0}(x+1)}{w''_{s,0}(x)},
  \label{Bd1..ds}\\
  D_{d_1\ldots d_s}(x)&\eqdef\alpha D'(x)
  \frac{w_s(x+1)}{w_s(x)}\frac{w''_{s,0}(x-1)}{w''_{s,0}(x)}.
  \label{Dd1..ds}
\end{align}
We have $B_{d_1\ldots d_s}(x)>0$ ($x\in\mathbb{Z}_{\geq 0}$),
$D_{d_1\ldots d_s}(x)>0$ ($x\in\mathbb{Z}_{\geq 1}$)
and $D_{d_1\ldots d_s}(0)=0$.
By using \eqref{I:id2} with $n=0$, we can show the relations
\begin{align*}
  B_{d_1\ldots d_s}(x)D_{d_1\ldots d_s}(x+1)
  &=\hat{B}_{d_1\ldots d_s}(x+1)\hat{D}_{d_1\ldots d_s}(x+1),\\
  B_{d_1\ldots d_s}(x)+D_{d_1\ldots d_s}(x)
  &=\hat{B}_{d_1\ldots d_s}(x)+\hat{D}_{d_1\ldots d_s}(x+1)
  +\tilde{\mathcal{E}}_{d_s}.
\end{align*}

It should be stressed that the above results after $s$-deletions are
independent of the orders of deletions ($\phi_{d_1\ldots d_s\,n}(x)$
and $\tilde{\phi}_{d_1\ldots d_s\,\text{v}}(x)$ may change sign).

\section{Multi-indexed Orthogonal Polynomials}
\label{sec:miop}

In this section we apply the method of virtual state deletions to
the exactly solvable systems whose eigenstates are described by
the Meixner, little $q$-Jacobi, and little $q$-Laguerre polynomials.
We delete $M$ virtual state vectors labeled by
\begin{equation}
  \mathcal{D}=\{d_1,d_2,\ldots,d_M\}
  \ \ (d_j\in\mathcal{V} : \text{mutually distinct}),
\end{equation}
and denote $\mathcal{H}_{d_1\ldots d_M}$, $\phi_{d_1\ldots d_M\,n}$,
$\mathcal{A}_{d_1\ldots d_M}$, etc. by $\mathcal{H}_{\mathcal{D}}$,
$\phi_{\mathcal{D}\,n}$, $\mathcal{A}_{\mathcal{D}}$, etc.
The Schr\"odinger equation for the multi-indexed Meixner, little $q$-Jacobi,
and little $q$-Laguerre systems reads
\begin{equation}
  \mathcal{H}_{\mathcal{D}}(\bm{\lambda})
  \phi_{\mathcal{D}\,n}(x;\bm{\lambda})
  =\mathcal{E}_n(\bm{\lambda})\phi_{\mathcal{D}\,n}(x;\bm{\lambda})
  \ \ (n=0,1,\ldots).
  \label{Mindeq}
\end{equation}
Without loss of generality, we assume $1\leq d_1<d_2<\cdots<d_M$.
As with the ordinary orthogonal polynomials, the index $n$ gives the
{\em number of nodes\/} of $\phi_{\mathcal{D}\,n}(x;\bm{\lambda})$
in $x\in(0,\infty)$.

\medskip

Let us denote the eigenvector $\phi_{\mathcal{D}\,n}(x;\bm{\lambda})$
in \eqref{phid1..dsn} after $M$-deletions ($s=M$) by
$\phi^{\text{gen}}_{\mathcal{D}\,n}(x;\bm{\lambda})$.
We define two polynomials $\check{\Xi}_{\mathcal{D}}(x;\bm{\lambda})$ and
$\check{P}_{\mathcal{D},n}(x;\bm{\lambda})$, to be called the denominator
polynomial and the {\em multi-indexed orthogonal polynomial\/}, respectively,
from the Casoratians as follows:
\begin{align}
  &\text{W}_{\text{C}}[\check{\xi}_{d_1},\ldots,\check{\xi}_{d_M}]
  (x;\bm{\lambda})
  =\mathcal{C}_{\mathcal{D}}(\bm{\lambda})\varphi_M(x)
  \check{\Xi}_{\mathcal{D}}(x;\bm{\lambda}),
  \label{XiDdef}\\
  &\text{W}_{\text{C}}[\check{\xi}_{d_1},\ldots,\check{\xi}_{d_M},
  \nu\check{P}_n](x;\bm{\lambda})
  =\mathcal{C}_{\mathcal{D},n}(\bm{\lambda})
  \varphi_{M+1}(x)
  \check{P}_{\mathcal{D},n}(x;\bm{\lambda})
  \nu(x;\bm{\lambda}+M\tilde{\bm{\delta}}),
  \label{cPDndef}\\
  &\tilde{\bm{\delta}}\eqdef\left\{
  \begin{array}{ll}
  (1,0)&:\text{M}\\
  (-1,1)&:\text{l$q$J}\\
  -1&:\text{l$q$L}
  \end{array}\right.,\ \quad
  \mathfrak{t}(\bm{\lambda})+u\bm{\delta}=
  \mathfrak{t}(\bm{\lambda}+u\tilde{\bm{\delta}})
  \ \ (\forall u\in\mathbb{R}).
  \label{deltat}
\end{align}
The constants $\mathcal{C}_{\mathcal{D}}(\bm{\lambda})$ and
$\mathcal{C}_{\mathcal{D},n}(\bm{\lambda})$ are specified later.
The auxiliary function $\varphi_M(x)$ is defined by \cite{os22}:
\begin{align}
  \varphi_M(x)&\eqdef\prod_{1\leq j<k\leq M}
  \frac{\eta(x+k-1)-\eta(x+j-1)}{\eta(k-j)}\qquad
  \bigl(\varphi_0(x)=\varphi_1(x)=1\bigr)\n
  &=\prod_{1\leq j<k\leq M}\varphi(x+j-1)
  =\left\{
  \begin{array}{ll}
  1&:\text{M}\\
  q^{\frac12M(M-1)x+\frac16M(M-1)(M-2)}&:\text{l$q$J,\,l$q$L}
  \end{array}\right..
  \label{varphiMdef}
\end{align}
The eigenvector \eqref{phid1..dsn} is rewritten as
\begin{align}
  \phi^{\text{gen}}_{\mathcal{D}\,n}(x;\bm{\lambda})
  &=(-1)^M\kappa^{\frac14M(M-1)}
  \frac{\mathcal{C}_{\mathcal{D},n}(\bm{\lambda})}
  {\mathcal{C}_{\mathcal{D}}(\bm{\lambda})}
  \sqrt{\prod_{j=1}^M\alpha(\bm{\lambda})
  B'\bigl(0;\bm{\lambda}+(j-1)\tilde{\bm{\delta}}\bigr)}\n
  &\quad\times
  \frac{\phi_0(x;\bm{\lambda}+M\tilde{\bm{\delta}})}
  {\sqrt{\check{\Xi}_{\mathcal{D}}(x;\bm{\lambda})
  \check{\Xi}_{\mathcal{D}}(x+1;\bm{\lambda})}}
  \check{P}_{\mathcal{D},n}(x;\bm{\lambda}).
  \label{phigen}
\end{align}
The denominator polynomial $\check{\Xi}_{\mathcal{D}}(x;\bm{\lambda})$
\eqref{XiDdef} and the multi-indexed orthogonal polynomial
$\check{P}_{\mathcal{D},n}(x;\bm{\lambda})$ \eqref{cPDndef} are polynomials
in $\eta$
\begin{equation}
  \check{\Xi}_{\mathcal{D}}(x;\bm{\lambda})\eqdef
  \Xi_{\mathcal{D}}\bigl(\eta(x);\bm{\lambda}\bigr),\quad
  \check{P}_{\mathcal{D},n}(x;\bm{\lambda})\eqdef
  P_{\mathcal{D},n}\bigl(\eta(x);\bm{\lambda}\bigr),
  \label{XiP_poly}
\end{equation}
and their {\em degrees\/} are generically $\ell_{\mathcal{D}}$ and
$\ell_{\mathcal{D}}+n$, respectively.
Here $\ell_{\mathcal{D}}$ is
\begin{equation}
  \ell_{\mathcal{D}}\eqdef\sum_{j=1}^Md_j-\tfrac12M(M-1).
  \label{lform}
\end{equation}
We adopt the universal normalisation for $\check{\Xi}_{\mathcal{D}}$ and
$\check{P}_{\mathcal{D},n}$:
$\check{\Xi}_{\mathcal{D}}(0;\bm{\lambda})=1$,
$\check{P}_{\mathcal{D},n}(0;\bm{\lambda})=1$,
which determine the constants $\mathcal{C}_{\mathcal{D}}(\bm{\lambda})$
and $\mathcal{C}_{\mathcal{D},n}(\bm{\lambda})$
(convention: $\prod\limits_{1\leq j<k\leq M}\!\!\!\!\!*=1$ for $M=1$),
\begin{align}
  \mathcal{C}_{\mathcal{D}}(\bm{\lambda})&\eqdef
  \frac{1}{\varphi_M(0)}
  \prod_{1\leq j<k\leq M}
  \frac{\tilde{\mathcal{E}}_{d_j}(\bm{\lambda})
  -\tilde{\mathcal{E}}_{d_k}(\bm{\lambda})}
  {\alpha(\bm{\lambda})B'(j-1;\bm{\lambda})},
  \label{CD}\\
  \mathcal{C}_{\mathcal{D},n}(\bm{\lambda})&\eqdef
  (-1)^M\mathcal{C}_{\mathcal{D}}(\bm{\lambda})
  \tilde{d}_{\mathcal{D},n}(\bm{\lambda})^2,\quad
  \tilde{d}_{\mathcal{D},n}(\bm{\lambda})^2\eqdef
  \frac{\varphi_M(0)}{\varphi_{M+1}(0)}
  \prod_{j=1}^M\frac{\mathcal{E}_n(\bm{\lambda})
  -\tilde{\mathcal{E}}_{d_j}(\bm{\lambda})}
  {\alpha(\bm{\lambda})B'(j-1;\bm{\lambda})}.
  \label{CDn}
\end{align}
The denominator polynomial $\check{\Xi}_{\mathcal{D}}(x;\bm{\lambda})$
is positive for $x\in\mathbb{Z}_{\geq 0}$.
The lowest degree multi-indexed orthogonal polynomial
$\check{P}_{\mathcal{D},0}(x;\bm{\lambda})$ is related to
$\check{\Xi}_{\mathcal{D}}(x;\bm{\lambda})$ by the parameter shift
$\bm{\lambda}\to\bm{\lambda}+\bm{\delta}$:
\begin{equation}
  \check{P}_{\mathcal{D},0}(x;\bm{\lambda})
  =\check{\Xi}_{\mathcal{D}}(x;\bm{\lambda}+\bm{\delta}).
  \label{PD0=XiD}
\end{equation}
The coefficients of the highest degree term of the polynomials $\check{P}_n$,
$\check{\Xi}_{\mathcal{D}}$ and $\check{P}_{\mathcal{D},n}$,
\begin{align}
  \check{P}_n(x;\bm{\lambda})&=c_n(\bm{\lambda})\eta(x)^n
  +\text{lower degree terms},\n
  \check{\Xi}_{\mathcal{D}}(x;\bm{\lambda})
  &=c^{\Xi}_{\mathcal{D}}(\bm{\lambda})\eta(x)^{\ell_{\mathcal{D}}}
  +\text{lower degree terms},\n
  \check{P}_{\mathcal{D},n}(x;\bm{\lambda})
  &=c^P_{\mathcal{D},n}(\bm{\lambda})\eta(x)^{\ell_{\mathcal{D}}+n}
  +\text{lower degree terms},
\end{align}
are
\begin{align}
  c_n(\bm{\lambda})&=\left\{
  \begin{array}{ll}
  (1-c^{-1})^n(\beta)_n^{-1}&:\text{M}\\[2pt]
  (-a)^{-n}q^{-n^2}(abq^{n+1};q)_n(bq;q)_n^{-1}&:\text{l$q$J}\\[2pt]
  (-a)^{-n}q^{-n^2}&:\text{l$q$L}
  \end{array}\right.,\\
  c^{\Xi}_{\mathcal{D}}(\bm{\lambda})
  &=\prod_{j=1}^M\frac{c_{d_j}\bigl(\mathfrak{t}(\bm{\lambda})\bigr)}
  {c_{j-1}\bigl(\mathfrak{t}(\bm{\lambda})\bigr)}
  \times\left\{
  \begin{array}{ll}
  1&:\text{M}\\
  {\displaystyle\prod_{1\leq j<k\leq M}\frac{aq^{-(j-1+k-1)}-bq}
  {aq^{-(d_j+d_k)}-bq}}&:\text{l$q$J}\\[14pt]
  {\displaystyle\prod_{1\leq j<k\leq M}q^{d_j+d_k-(j-1+k-1)}}&:\text{l$q$L}
  \end{array}\right.,\\
  c^P_{\mathcal{D},n}(\bm{\lambda})
  &=c^{\Xi}_{\mathcal{D}}(\bm{\lambda})c_n(\bm{\lambda})
  \times\left\{
  \begin{array}{ll}
  {\displaystyle\prod_{j=1}^M\frac{\beta+j-1}{\beta+d_j+n}}&:\text{M}\\
  {\displaystyle\prod_{j=1}^M\frac{q^{-(j-1)}-bq}{q^{-(d_j+n)}-bq}}
  &:\text{l$q$J}\\[14pt]
  {\displaystyle\prod_{j=1}^Mq^{d_j+n-(j-1)}}&:\text{l$q$L}
  \end{array}\right..
\end{align} 
Among these coefficients, $c^{\Xi}_{\mathcal{D}}(\bm{\lambda})$ and
$c^P_{\mathcal{D},n}(\bm{\lambda})$ for the l$q$J can vanish for certain
ratios of the parameters, that is $a=bq^{m+1}$, $m\in\mathbb{Z}_{\ge0}$ as
remarked below \eqref{lqJvirt}.
These special configurations of parameters require separate treatment.

The potential functions $B_{\mathcal{D}}(x)$ and $D_{\mathcal{D}}(x)$
\eqref{Bd1..ds}--\eqref{Dd1..ds} after $M$-deletion ($s=M$) can be
expressed neatly in terms of the denominator polynomial:
\begin{align}
  B_{\mathcal{D}}(x;\bm{\lambda})&=B(x;\bm{\lambda}+M\tilde{\bm{\delta}})\,
  \frac{\check{\Xi}_{\mathcal{D}}(x;\bm{\lambda})}
  {\check{\Xi}_{\mathcal{D}}(x+1;\bm{\lambda})}
  \frac{\check{\Xi}_{\mathcal{D}}(x+1;\bm{\lambda}+\bm{\delta})}
  {\check{\Xi}_{\mathcal{D}}(x;\bm{\lambda}+\bm{\delta})},
  \label{BD2}\\
  D_{\mathcal{D}}(x;\bm{\lambda})&=D(x)\,
  \frac{\check{\Xi}_{\mathcal{D}}(x+1;\bm{\lambda})}
  {\check{\Xi}_{\mathcal{D}}(x;\bm{\lambda})}
  \frac{\check{\Xi}_{\mathcal{D}}(x-1;\bm{\lambda}+\bm{\delta})}
  {\check{\Xi}_{\mathcal{D}}(x;\bm{\lambda}+\bm{\delta})},
  \label{DD2}\\
  \mathcal{A}_{\mathcal{D}}(\bm{\lambda})&=
  \sqrt{B_{\mathcal{D}}(x;\bm{\lambda})}
  -e^{\partial}\sqrt{D_{\mathcal{D}}(x;\bm{\lambda})},\quad
  \mathcal{A}_{\mathcal{D}}^{\dagger}(\bm{\lambda})
  =\sqrt{B_{\mathcal{D}}(x;\bm{\lambda})}
  -\sqrt{D_{\mathcal{D}}(x;\bm{\lambda})}\,e^{-\partial}.
\end{align}
The ground state eigenvector $\phi_{\mathcal{D}\,0}(x)$ is expressed by
$\phi_0(x)$ \eqref{phi0} and
$\check{\Xi}_{\mathcal{D}}(x)$:
\begin{align}
  \phi_{\mathcal{D}\,0}(x;\bm{\lambda})&=
  \prod_{y=0}^{x-1}\sqrt{\frac{B_{\mathcal{D}}(y;\bm{\lambda})}
  {D_{\mathcal{D}}(y+1;\bm{\lambda})}}
  =\phi_0(x;\bm{\lambda}+M\tilde{\bm{\delta}})
  \sqrt{\frac{\check{\Xi}_{\mathcal{D}}(1;\bm{\lambda})}
  {\check{\Xi}_{\mathcal{D}}(x;\bm{\lambda})
  \check{\Xi}_{\mathcal{D}}(x+1;\bm{\lambda})}}\,
  \check{\Xi}_{\mathcal{D}}(x;\bm{\lambda}+\bm{\delta})\n
  &=\psi_{\mathcal{D}}(x;\bm{\lambda})
  \check{P}_{\mathcal{D},0}(x;\bm{\lambda})
  \propto\phi^{\text{gen}}_{\mathcal{D}\,0}(x;\bm{\lambda}),\\
  \psi_{\mathcal{D}}(x;\bm{\lambda})&\eqdef
  \sqrt{\check{\Xi}_{\mathcal{D}}(1;\bm{\lambda})}\,
  \frac{\phi_0(x;\bm{\lambda}+M\tilde{\bm{\delta}})}
  {\sqrt{\check{\Xi}_{\mathcal{D}}(x;\bm{\lambda})\,
  \check{\Xi}_{\mathcal{D}}(x+1;\bm{\lambda})}},\quad
  \psi_{\mathcal{D}}(0;\bm{\lambda})=1.
\end{align}
We arrive at the normalised eigenvector
$\phi_{\mathcal{D}\,n}(x;\bm{\lambda})$ with the orthogonality relation,
\begin{align}
  &\phi_{\mathcal{D}\,n}(x;\bm{\lambda})
  \eqdef\psi_{\mathcal{D}}(x;\bm{\lambda})
  \check{P}_{\mathcal{D},n}(x;\bm{\lambda})
  \propto\phi^{\text{gen}}_{\mathcal{D}\,n}(x;\bm{\lambda}),\quad
  \phi_{\mathcal{D}\,n}(0;\bm{\lambda})=1,\\
  &\sum_{x=0}^{\infty}
  \frac{\psi_{\mathcal{D}}(x;\bm{\lambda})^2}
  {\check{\Xi}_{\mathcal{D}}(1;\bm{\lambda})}
  \check{P}_{\mathcal{D},n}(x;\bm{\lambda})
  \check{P}_{\mathcal{D},m}(x;\bm{\lambda})
  =\frac{\delta_{nm}}{d_n(\bm{\lambda})^2
  \tilde{d}_{\mathcal{D},n}(\bm{\lambda})^2}
  \ \ (n,m\in\mathbb{Z}_{\geq 0}).
\end{align}
It is worthwhile to emphasise that the above orthogonality relation is
a rational equation of $\bm{\lambda}$ or $q^{\bm{\lambda}}$, and it is
valid for any value of $\bm{\lambda}$ (except for the zeros of denominators
and those parameter ranges giving divergent sums)
but the weight function may not be positive definite.

The shape invariance of the original system is inherited by the deformed
systems:
\begin{equation}
  \mathcal{A}_{\mathcal{D}}(\bm{\lambda})
  \mathcal{A}_{\mathcal{D}}(\bm{\lambda})^{\dagger}
  =\kappa\mathcal{A}_{\mathcal{D}}(\bm{\lambda}+\bm{\delta})^{\dagger}
  \mathcal{A}_{\mathcal{D}}(\bm{\lambda}+\bm{\delta})
  +\mathcal{E}_1(\bm{\lambda}).
  \label{shapeinvD}
\end{equation}
As a consequence of the shape invariance and the normalisation,
the actions of $\mathcal{A}_{\mathcal{D}}(\bm{\lambda})$ and
$\mathcal{A}_{\mathcal{D}}(\bm{\lambda})^{\dagger}$ on the eigenvectors
$\phi_{\mathcal{D}\,n}(x;\bm{\lambda})$ are
\begin{align}
  &\mathcal{A}_{\mathcal{D}}(\bm{\lambda})
  \phi_{\mathcal{D}\,n}(x;\bm{\lambda})
  =\frac{\mathcal{E}_n(\bm{\lambda})}
  {\sqrt{B_{\mathcal{D}}(0;\bm{\lambda})}}\,
  \phi_{\mathcal{D}\,n-1}(x;\bm{\lambda}+\bm{\delta}),
  \label{ADphiDn=}\\
  &\mathcal{A}_{\mathcal{D}}(\bm{\lambda})^{\dagger}
  \phi_{\mathcal{D}\,n-1}(x;\bm{\lambda}+\bm{\delta})
  =\sqrt{B_{\mathcal{D}}(0;\bm{\lambda})}\,
  \phi_{\mathcal{D}\,n}(x;\bm{\lambda}).
  \label{ADdphiDn=}
\end{align}
The forward and backward shift operators are defined by
\begin{align}
  \mathcal{F}_{\mathcal{D}}(\bm{\lambda})&\eqdef
  \sqrt{B_{\mathcal{D}}(0;\bm{\lambda})}\,
  \psi_{\mathcal{D}}(x;\bm{\lambda}+\bm{\delta})^{-1}\circ
  \mathcal{A}_{\mathcal{D}}(\bm{\lambda})\circ
  \psi_{\mathcal{D}}(x;\bm{\lambda})\n
  &=\frac{B(0;\bm{\lambda}+M\tilde{\bm{\delta}})}
  {\varphi(x)\check{\Xi}_{\mathcal{D}}(x+1;\bm{\lambda})}
  \Bigl(\check{\Xi}_{\mathcal{D}}(x+1;\bm{\lambda}+\bm{\delta})
  -\check{\Xi}_{\mathcal{D}}(x;\bm{\lambda}+\bm{\delta})e^{\partial}\Bigr),
  \label{calFD}\\
  \mathcal{B}_{\mathcal{D}}(\bm{\lambda})&\eqdef
  \frac{1}{\sqrt{B_{\mathcal{D}}(0;\bm{\lambda})}}\,
  \psi_{\mathcal{D}}(x;\bm{\lambda})^{-1}\circ
  \mathcal{A}_{\mathcal{D}}(\bm{\lambda})^{\dagger}\circ
  \psi_{\mathcal{D}}(x;\bm{\lambda}+\bm{\delta})\n
  &=\frac{1}{B(0;\bm{\lambda}+M\tilde{\bm{\delta}})
  \check{\Xi}_{\mathcal{D}}(x;\bm{\lambda}+\bm{\delta})}
  \Bigl(B(x;\bm{\lambda}+M\tilde{\bm{\delta}})
  \check{\Xi}_{\mathcal{D}}(x;\bm{\lambda})
  -D(x)\check{\Xi}_{\mathcal{D}}(x+1;\bm{\lambda})e^{-\partial}\Bigr)
  \varphi(x)\n
  &=\psi_{\mathcal{D}}(x;\bm{\lambda})^{-2}\circ
  \frac{\check{\Xi}_{\mathcal{D}}(1;\bm{\lambda})}
  {\check{\Xi}_{\mathcal{D}}(1;\bm{\lambda}+\bm{\delta})}
  \Bigl(\check{\Xi}_{\mathcal{D}}(x+1;\bm{\lambda}+\bm{\delta})
  -\check{\Xi}_{\mathcal{D}}(x-1;\bm{\lambda}+\bm{\delta})e^{-\partial}
  \Bigr)\n
  &\qquad\qquad\qquad\qquad\times
  \frac{1}{\varphi(x)\check{\Xi}_{\mathcal{D}}(x+1;\bm{\lambda})}\circ
  \psi_{\mathcal{D}}(x;\bm{\lambda}+\bm{\delta})^2,
  \label{calBD}
\end{align}
and their actions on $\check{P}_{\mathcal{D},n}(x;\bm{\lambda})$ are
\begin{equation}
  \mathcal{F}_{\mathcal{D}}(\bm{\lambda})
  \check{P}_{\mathcal{D},n}(x;\bm{\lambda})
  =\mathcal{E}_n(\bm{\lambda})
  \check{P}_{\mathcal{D},n-1}(x;\bm{\lambda}+\bm{\delta}),\quad
  \mathcal{B}_{\mathcal{D}}(\bm{\lambda})
  \check{P}_{\mathcal{D},n-1}(x;\bm{\lambda}+\bm{\delta})
  =\check{P}_{\mathcal{D},n}(x;\bm{\lambda}).
  \label{BDPDn=}
\end{equation}
The similarity transformed Hamiltonian is
\begin{align}
  \widetilde{\mathcal{H}}_{\mathcal{D}}(\bm{\lambda})
  &\eqdef\psi_{\mathcal{D}}(x;\bm{\lambda})^{-1}\circ
  \mathcal{H}_{\mathcal{D}}(\bm{\lambda})\circ
  \psi_{\mathcal{D}}(x;\bm{\lambda})
  =\mathcal{B}_{\mathcal{D}}(\bm{\lambda})
  \mathcal{F}_{\mathcal{D}}(\bm{\lambda})\n
  &=B(x;\bm{\lambda}+M\tilde{\bm{\delta}})\,
  \frac{\check{\Xi}_{\mathcal{D}}(x;\bm{\lambda})}
  {\check{\Xi}_{\mathcal{D}}(x+1;\bm{\lambda})}
  \biggl(\frac{\check{\Xi}_{\mathcal{D}}(x+1;\bm{\lambda}+\bm{\delta})}
  {\check{\Xi}_{\mathcal{D}}(x;\bm{\lambda}+\bm{\delta})}-e^{\partial}
  \biggr)\n
  &\quad+D(x)\,
  \frac{\check{\Xi}_{\mathcal{D}}(x+1;\bm{\lambda})}
  {\check{\Xi}_{\mathcal{D}}(x;\bm{\lambda})}
  \biggl(\frac{\check{\Xi}_{\mathcal{D}}(x-1;\bm{\lambda}+\bm{\delta})}
  {\check{\Xi}_{\mathcal{D}}(x;\bm{\lambda}+\bm{\delta})}-e^{-\partial}
  \biggr),
  \label{MHtdef}
\end{align}
and the multi-indexed orthogonal polynomials
$\check{P}_{\mathcal{D},n}(x;\bm{\lambda})$ are its eigenpolynomials:
\begin{equation}
  \widetilde{\mathcal{H}}_{\mathcal{D}}(\bm{\lambda})
  \check{P}_{\mathcal{D},n}(x;\bm{\lambda})=\mathcal{E}_n(\bm{\lambda})
  \check{P}_{\mathcal{D},n}(x;\bm{\lambda}).
  \label{tHPDn=}
\end{equation}

Including the level 0 deletion corresponds to $M-1$ virtual
states deletion:
\begin{equation}
  \check{P}_{\mathcal{D},n}(x;\bm{\lambda})\Bigm|_{d_M=0}
  =\check{P}_{\mathcal{D}',n}(x;\bm{\lambda}+\tilde{\bm{\delta}}),\quad
  \mathcal{D}'=\{d_1-1,\ldots,d_{M-1}-1\}.
  \label{dM=0}
\end{equation}
The denominator polynomial $\Xi_{\mathcal{D}}$ behaves similarly.
This is why we have restricted $d_j\geq 1$.

\subsection{Limits}
\label{sec:limit}

It is well known \cite{kls} that the Meixner \eqref{MeixnerP},
little $q$-Jacobi \eqref{littleqjacobinorm} and little $q$-Laguerre
\eqref{littleqlaguerrenorm} polynomials reduce to the Laguerre (L)
$L^{(\alpha)}_n(\eta)$ or Jacobi (J) $P^{(\alpha,\beta)}_n(\eta)$ polynomials
in certain limits:
\begin{alignat}{2}
  \text{M}:&\ \ \lim_{c\to 1}P_n\Bigl(\frac{\eta}{1-c};(\alpha+1,c)\Bigr)
  =\frac{L_n^{(\alpha)}(\eta)}{L_n^{(\alpha)}(0)},
  &&\text{\em i.e.}\ \ \beta=\alpha+1,
  \label{M->L}\\
  \text{l$q$J}:&\ \ \lim_{q\to 1}P_n\bigl(1-\eta;(\alpha,\beta)\bigr)
  =\frac{P^{(\alpha,\beta)}_n(1-2\eta)}{P^{(\alpha,\beta)}_n(-1)},
  &\quad\ \ &\text{\em i.e.}\ \ a=q^{\alpha},\ b=q^{\beta},
  \label{lqJ->J}\\
  \text{l$q$L}:&\ \ \lim_{q\to 1}
  \frac{P_n\bigl(1-(1-q)\eta;\alpha\bigr)}{P_n(1;\alpha)}
  =\frac{L^{(\alpha)}_n(\eta)}{L^{(\alpha)}_n(0)},
  &&\text{\em i.e.}\ \ a=q^{\alpha}.
  \label{lqL->L}
\end{alignat}
The constant factors on r.h.s.\ reflect the universal normalisation
$P_n(0)=1$ of M, l$q$J(L) polynomials \eqref{Pzero}.
Note that $L^{(\alpha)}_n(0)=\frac{1}{n!}(\alpha+1)_n$ and
$P^{(\alpha,\beta)}_n(-1)=(-1)^nP^{(\beta,\alpha)}_n(1)
=\frac{(-1)^n}{n!}(\beta+1)_n$.

These relations \eqref{M->L}--\eqref{lqL->L} also mean that the deforming
polynomials $\xi(\eta;\bm{\lambda})$ of multi-indexed M, l$q$J(L) polynomials,
\eqref{Mvirt}, \eqref{lqJvirt} reduce to those of the multi-indexed Laguerre
and Jacobi polynomials:
\begin{align}
  \text{M}:&\ \ \lim_{c\to 1}\xi_{\text{v}}\Bigl(\frac{\eta}{1-c};
  (\alpha+1,c)\Bigr)
  =\frac{L_{\text{v}}^{(\alpha)}(-\eta)}{L_{\text{v}}^{(\alpha)}(0)}, 
  \label{xiM->L}\\
  \text{l$q$J}:&\ \ \lim_{q\to 1}\xi_{\text{v}}\bigl(1-\eta;
  (\alpha,\beta)\bigr)
  =\frac{P^{(-\alpha,\beta)}_{\text{v}}(1-2\eta)}
  {P^{(-\alpha,\beta)}_{\text{v}}(-1)},
  \label{xilqJ->J}\\
  \text{l$q$L}:&\ \ \lim_{q\to 1}
  \frac{\xi_{\text{v}}\bigl(1-(1-q)\eta;\alpha\bigr)}{\xi_{\text{v}}(1;\alpha)}
  =\frac{L^{(-\alpha)}_{\text{v}}(\eta)}{L^{(-\alpha)}_{\text{v}}(0)}.
  \label{xilqL->L}
\end{align}
As shown in \cite{os25}, the deforming polynomial
$L_{\text{v}}^{(\alpha)}(-\eta)$ \eqref{xiM->L} will lead to type $\I$
multi-indexed Laguerre polynomials, whereas $L^{(-\alpha)}_{\text{v}}(\eta)$
\eqref{xilqL->L} to type $\II$ Laguerre.
The deforming polynomial $P^{(-\alpha,\beta)}_{\text{v}}(\eta)$ will
generate type $\II$ multi-indexed Jacobi polynomials.
Based on these formulas \eqref{M->L}--\eqref{xilqL->L} and the Casoratian
definition \eqref{cPDndef}, it is straightforward to demonstrate that the
multi-indexed M, l$q$J(L) polynomials reduce to the Wronskian definition of
the corresponding multi-indexed Laguerre and Jacobi polynomials \cite{os25}:
\begin{alignat}{2}
  \text{M}:&\ \ \lim_{c\to 1}P_{\mathcal{D},n}
  \Bigl(\frac{\eta}{1-c};(\alpha+1,c)\Bigr)
  =\frac{P^{\text{L}}_{\mathcal{D},n}(\eta;g)}
  {P^{\text{L}}_{\mathcal{D},n}(0;g)},
  &&g\eqdef\alpha+\tfrac12,
  \label{D,M->L}\\
  \text{l$q$J}:&\ \ \lim_{q\to 1}P_{\mathcal{D},n}
  \bigl(1-\eta;(\alpha,\beta)\bigr)
  =\frac{P^{\text{J}}_{\mathcal{D},n}\bigl(1-2\eta;(g,h)\bigr)}
  {P^{\text{J}}_{\mathcal{D},n}\bigl(-1;(g,h)\bigr)},
  &\quad\ \ &g\eqdef\alpha+\tfrac12,\ h\eqdef\beta+\tfrac12,
  \label{D,lqJ->J}\\
  \text{l$q$L}:&\ \ \lim_{q\to 1}
  \frac{P_{\mathcal{D},n}\bigl(1-(1-q)\eta;\alpha\bigr)}
  {P_{\mathcal{D},n}\bigl(1;\alpha\bigr)}
  =\frac{P^{\text{L}}_{\mathcal{D},n}(\eta;g)}
  {P^{\text{L}}_{\mathcal{D},n}(0;g)},
  &&g\eqdef\alpha+\tfrac12.
  \label{D,lqL->L}
\end{alignat}
Here $P^{\text{L}}_{\mathcal{D},n}(\eta;\bm{\lambda})$ and
$P^{\text{J}}_{\mathcal{D},n}(\eta;\bm{\lambda})$ are the multi-indexed
Laguerre and Jacobi polynomials \cite{os25} with the index set
$\mathcal{D}=\{d_1,\ldots,d_M\}$, which are all type $\I$ for \eqref{D,M->L}
and all type $\II$ for \eqref{D,lqJ->J}--\eqref{D,lqL->L}.
Again the constant factors on r.h.s.\ reflect the universal normalisation 
of the multi-indexed M, l$q$J(L) polynomials.

\section{Summary and Comments}
\label{sec:summary}

Following the general procedure \cite{os25}, {\em i.e.\/}
\romannumeral1) generating virtual state vectors by twisting parameters,
\romannumeral2) applying multiple Darboux transformations by using the above
virtual state vectors as seed solutions, the multi-indexed Meixner, little
$q$-Jacobi and little $q$-Laguerre polynomials are constructed.
As with the other types of multi-indexed orthogonal polynomials, they are
shape invariant and form a complete basis in the $\ell^2$ Hilbert space,
although they lack certain lower degrees.

For the quantum mechanical systems described by the multi-index Laguerre
and Jacobi polynomials, the eigenfunctions have the following form \cite{os25}
\begin{equation*}
  \phi_{\mathcal{D}\,n}(x;\bm{\lambda})
  \propto\frac{\phi_0(x;\bm{\lambda}')}
  {\check{\Xi}_{\mathcal{D}}(x;\bm{\lambda})}
  \check{P}_{\mathcal{D},n}(x;\bm{\lambda})
  \ \biggl(=\phi_0(x;\bm{\lambda}')
  \frac{\check{P}_{\mathcal{D},n}(x;\bm{\lambda})}
  {\check{\Xi}_{\mathcal{D}}(x;\bm{\lambda})}\biggr),
\end{equation*}
where $\bm{\lambda}'=\bm{\lambda}+(\text{shift})$.
Usual interpretation is that the polynomials
$\check{P}_{\mathcal{D},n}(x;\bm{\lambda})
=P_{\mathcal{D},n}(\eta(x);\bm{\lambda})$
are orthogonal with respect to the weight function
$\bigl(\phi_0(x;\bm{\lambda}')/
\check{\Xi}_{\mathcal{D}}(x;\bm{\lambda})\bigr)^2$ and form a complete set.
Another interpretation is possible; the rational functions
$\check{P}_{\mathcal{D},n}(x;\bm{\lambda})/
\check{\Xi}_{\mathcal{D}}(x;\bm{\lambda})$
are orthogonal with respect to the original weight function (with shifted
parameters) $\bigl(\phi_0(x;\bm{\lambda}')\bigr)^2$ and form a complete set.
The lowest degree element is almost constant
$\check{P}_{\mathcal{D},0}(x;\bm{\lambda})/
\check{\Xi}_{\mathcal{D}}(x;\bm{\lambda})\propto
\check{\Xi}_{\mathcal{D}}(x;\bm{\lambda}+\bm{\delta})/
\check{\Xi}_{\mathcal{D}}(x;\bm{\lambda})$ as \eqref{PD0=XiD}.
The ($n+1$)-th element has $n$ zeros in $x\in(0,\infty)$
having the structure
$(\text{degree $\ell_{\mathcal{D}}+n$ polynomial})/
(\text{degree $\ell_{\mathcal{D}}$ polynomial})$.
For the discrete quantum mechanics with real shifts described by the
multi-index polynomials, e.g.\ ($q$-)Racah \cite{os26} and M, l$q$J and l$q$L
in this paper, the above rational functions are replaced as \eqref{phigen}
\begin{equation*}
  \frac{\check{P}_{\mathcal{D},n}(x;\bm{\lambda})}
  {\check{\Xi}_{\mathcal{D}}(x;\bm{\lambda})}
  \ \to\ \frac{\check{P}_{\mathcal{D},n}(x;\bm{\lambda})}
  {\sqrt{\check{\Xi}_{\mathcal{D}}(x;\bm{\lambda})
  \check{\Xi}_{\mathcal{D}}(x+1;\bm{\lambda})}}.
\end{equation*}
These functions $\check{P}_{\mathcal{D},n}(x;\bm{\lambda})/
\sqrt{\check{\Xi}_{\mathcal{D}}(x;\bm{\lambda})
\check{\Xi}_{\mathcal{D}}(x+1;\bm{\lambda})}$ are no longer rational functions
but they are orthogonal with respect to the original weight function
(with shifted parameters)
$\bigl(\phi_0(x;\bm{\lambda}+M\tilde{\bm{\delta}})\bigr)^2$
and form a complete set.

These multi-indexed orthogonal polynomials reduce to the multi-indexed
Laguerre or Jacobi polynomials in certain limits:
M $\to$ L, l$q$J $\to$ J, l$q$L $\to$ L.
It is an interesting question whether the other types of multi-indexed
orthogonal polynomials, {\em i.e.\/} type $\II$ for M and type $\I$ for
l$q$J(L), could be constructed for the M, l$q$J and l$q$L systems.
The existence of such multi-indexed orthogonal polynomials is not obvious
at all, as the discrete symmetries of the reduced systems do not imply
those of the original difference Schr\"odinger systems.

{}From exactly solvable models of discrete QM with real shifts, exactly
solvable birth and death processes can be constructed \cite{bdproc,os34}.
The exactly solvable models described by the multi-indexed M, l$q$J, l$q$L
polynomials provide new exactly solvable birth and death processes.

Dur\'an \cite{d14} reported multi-indexed Meixner polynomials, in which
he used both eigenvectors and virtual state vectors as seed solutions.
The method of construction consists in dualizing Krall discrete orthogonal
polynomials. Krall discrete orthogonal polynomials are sequences
$(p_n)_{n\in{\mathbb N}}$ of orthogonal polynomials which in addition are
eigenfunctions of a higher order difference operator. 
Dur\'an uses suitable instances of Christoffel transforms for constructing
Krall discrete orthogonal polynomials \cite{dn1,dn2,dn3}.
Using eigenvectors as seed solutions for constructing exactly solvable QM
and multi-indexed orthogonal polynomials was established by Krein \cite{krein}
and Adler \cite{adler} more than 20 years ago, and for discrete QM with
real shifts it was done in \cite{os22}.
Dur\'an suggested an interesting possibility that the parameter ranges
for positive definite orthogonality measure of the multi-indexed Meixner
polynomials could be enlarged when the eigenvectors and virtual state vectors
are used simultaneously.
In this paper we have not pursued its feasibility, as we have not adopted
eigenvectors as seed solutions.

There are 39 classical orthogonal polynomials satisfying second order
difference equations in the Askey scheme of hypergeometric orthogonal
polynomials and its $q$-version.
The aim of our multi-indexed orthogonal polynomial project is to construct
the {\em shape invariant deformation} of all possible examples in the
($q$-)Askey scheme.
Obviously, the continuous $q$-Hermite polynomial is excluded as it has no
twistable parameter, like the Hermite polynomial.
We have reported the shape invariant multi-indexed deformations of the
Wilson, Askey-Wilson \cite{os27}, Racah and $q$-Racah \cite{os26} polynomials
on top of those reported in the present paper.
All the rest of the polynomials in the ($q$-)Askey scheme are obtained by
parameter restriction and/or certain limiting procedures from the above
four polynomials.
For some of them, the multi-indexed deformation is expected to be compatible
with the parameter restriction/limiting procedures, but the actual
verification is yet to be done for each case.
Among these, the (dual) big $q$-Jacobi (Laguerre), (dual) Al-Salam-Carltz
$\I$ and the discrete $q$-Hermite polynomials provide a big challenge.
The orthogonality measures of these polynomials are of the Jackson integral
\cite{askey} type and their quantum mechanical formulation requires
two-component Hamiltonians \cite{os34}.
The parameter twisting is rather complicated in these cases.
It seems brand new thinking is required to tackle this problem.

\section*{Acknowledgements}

S.\,O. is supported in part by Grant-in-Aid for Scientific Research
from the Ministry of Education, Culture, Sports, Science and Technology
(MEXT), No.25400395.

\goodbreak

\end{document}